\documentclass[12pt,a4paper]{article}
\usepackage{times}
\usepackage{graphicx}
\usepackage[square,numbers,comma,compress]{natbib}
\usepackage[fleqn]{amsmath}
\usepackage{amssymb,amsfonts,amscd,stmaryrd,siunitx,empheq}
\usepackage{xcolor}
\usepackage{hyperref}
\hypersetup{colorlinks=true,linkcolor=red,citecolor=blue}
\newcommand{\bm}[1]{\boldsymbol{#1}}

\newcommand{\jump}[1]{\left\llbracket #1 \right\rrbracket}
\newcommand{\mean}[1]{\left\{ #1 \right\}}
\newcommand{\grad}{\bm{\nabla}}
\newcommand{\hsigma}{\widehat{\bm{\sigma}}}
\newcommand{\tsigma}{\widetilde{\bm{\sigma}}}
\newcommand{\bu}{\bm{u}}
\newcommand{\bv}{\bm{v}}
\newcommand{\bE}{\bm{E}}
\newcommand{\I}{\mathcal{I}}
\newcommand{\hD}{\widehat{\bm{D}}}
\newcommand{\nsd}{{n_\text{sd}}}
\newcommand{\numel}{{n_\text{el}}}

\newcommand{\nedges}{{n_\text{edg}}}
\newcommand{\norm}[1]{\left\lVert#1\right\rVert_{\mathcal{L}_2(\mathcal{I})}}

\newtheorem{remark}{Remark}

\title{A C0 interior penalty finite element method for flexoelectricity}
\author
{J. Ventura$^{1,a}$, D. Codony$^{1,b}$, S. Fern\'andez-M\'endez$^{1,c,\ast}$\\
	\\
	\small{
	$^1$ Laboratori de Càlcul Numèric, E.T.S. de Ingeniería de
Caminos,\vspace{-.4em}}
	\\
	\small{Universitat Politècnica de Catalunya, Barcelona,
Spain }
	\\\vspace{-.4em}
	\small{$^a$ Email: jordi.ventura.siches@estudiant.upc.edu}
	\\\vspace{-.4em}
	\small{$^b$ Email: david.codony@upc.edu}
	\\
	\small{$^c$ Email: sonia.fernandez@upc.edu}
	\\
	\small{$^\ast$ Corresponding author}
}
\date{}
\begin{document}
\maketitle
\newcommand{\sep}{,~}
\begin{abstract} 
We propose a $\mathcal{C}^0$ Interior Penalty Method (C0-IPM) for the computational modelling of flexoelectricity, with application also to strain gradient elasticity, as a simplified case. Standard high-order $\mathcal{C}^0$ finite element approximations, with nodal basis, are considered. The proposed C0-IPM formulation involves second derivatives in the interior of the elements, plus integrals on the mesh faces (sides in 2D), that impose $\mathcal{C}^1$ continuity of the displacement in weak form.  The formulation is stable for large enough interior penalty parameter, which can be estimated solving an eigenvalue problem. The applicability and convergence of the method is demonstrated with 2D and 3D numerical examples.
\end{abstract}

\section{Introduction}
\label{intro}

The rising interest on microtechnology  evidences the need for mathematical and computational models suitable for small scales, often giving rise to $4^\text{th}$ order Partial Differential Equations (PDEs). In particular flexoelectric effects become relevant, and may be crucial, in the design of small electromechanical devices or for the understanding of physical phenomena \cite{Shu2019}.
The modelling of flexoelectricity involves a two-way coupling between strain gradient and electric field; strain gradient elasticity is frequently also included in the model to regularise the problem, leading to a system of $4^\text{th}$ order PDEs. 

Several numerical strategies have been recently proposed for the solution of  flexoelectricity problems, based on the use of $\mathcal{C}^1$ approximation spaces or on mixed formulations.  Mixed formulations split the PDE in two $2^\text{nd}$ order PDEs, allowing the use of $\mathcal{C}^0$ Finite Element (FE) approximations \cite{Mao2016,Deng2017}. The approximation spaces, for the primal unknown and for the additional unknowns, must fulfill some conditions for stability that lead to approximation spaces with cumbersome definitions, and difficult extension to 3D or high-order approximations. However, the main drawback of mixed formulations in the high computational cost due to the additional unknowns. 

On other hand, $\mathcal{C}^1$ approximations can be directly used for the discretization of the weak form in $\mathcal{H}^2$, involving $2^\text{nd}$ order derivatives, without additional unknowns. The first successful attempt in this direction considered a meshless method: the maximum entropy method \cite{Abdollahi2014}. 
Unfortunately, the computational cost of accurate meshless methods is high, mainly due to the excessive number of integration points for accurate solutions, and the large stencils in the discrete matrices. Aiming to improve the efficiency, a solution based on Isogeometric Analysis (IGA) is proposed in \cite{Ghasemi2017}. The solution is approximated by means of Non-Uniform Rational B-Splines (NURBS). An interesting critical comparison of IGA, meshless and mixed methods for flexoelectricity can be found at \cite{Nanthakumar2017}, although numerical examples are restricted to simple geometries. The conclusion is that IGA is very efficient on regular grids, corresponding to a transformation of a rectangle grid, where plain B-Spline approximations can be easily defined. 
However, in a more general context, defining a NURBS approximation with $\mathcal{C}^1$ continuity in a whole domain with complex shape may not be straightforward, and the numerical integration of the resulting NURBS may be, again, very expensive \cite{Sevilla2011}. 
An efficient alternative for complex domains is the immersed B-Spline method proposed in \cite{Codony2019}. It considers B-Spline approximations based on a background regular grid, with an embedded domain. The applicability of the proposal is demonstrated with 2D and 3D complex geometries. The weak points of immersed B-Splines are the usual ones in the context of embedded domains: the robust definition of numerical integration in cropped cells (intersected by the domain boundary), which is specially challenging in 3D, and the ill-conditioning problems in the presence of cells with a small portion in the domain, that can be alleviated with specific techniques \cite{dePrenter2017}.

Here we propose a $\mathcal{C}^0$ Interior Penalty Method (C0-IPM) for the solution of flexoelectricity. A standard $\mathcal{C}^0$ FE approximation is considered, and $\mathcal{C}^1$ continuity {between elements} is imposed in weak form by means of the Interior Penalty Method (IPM). The procedure for the derivation of the C0-IPM weak form is analogous to the derivation of the IPM in the context of DG methods for $2^\text{nd}$ order PDEs \cite{Arnold1982}, or Nitsche's method for weak imposition of Dirichlet boundary conditions \cite{Nitsche1971}, but now applied to the continuity of normal derivatives on element boundaries.
The resulting weak form involves second derivatives in the interior of the elements, plus integrals on the faces (sides in 2D) of the mesh, that impose $\mathcal{C}^1$ continuity of the displacement in weak form. 
C0-IPM formulations overcome the disadvantages of other methods, because they allow the use of standard $\mathcal{C}^0$ FE approximations. Namely, {\it (i)} the computational mesh can be adapted to any geometry, with localised refinement were needed, {\it(ii)} there is no need to use embedded discretizations, avoiding the consequent ill-conditioning problems and the definition of special numerical integration for cropped elements, {\it(iii)} there are no additional unknowns and {\it(iv)} they handle material interfaces in a natural way. In summary, C0-IPM retains the computational efficiency and the versatility that make standard FEs the preferred method for many practitioners in the computational mechanics community.

In \cite{Engel2002}, C0-IPM formulations, there referred to as continuous/discontinuous finite elements, are applied to several problems modelled by $4^\text{th}$ order PDEs, including Kirchhoff plates and 1D strain gradient elasticity.  Numerical experiments show the applicability of the formulation in both applications, but convergence studies are limited to 1D examples. The C0-IPM formulation is then analysed {in \cite{Brenner2005}} for the 2D biharmonic equation, with first and second Dirichlet conditions, including a convergence analysis that shows that the method is convergent for $p\geq 2$, but may have suboptimal convergence depending on the degree and the penalty parameter. An experimental convergence study for Kirchhoff plates can be found at \cite{Fojo2020}. The numerical results demonstrate the applicability of the method for degree greater or equal to $3$, and also show slightly suboptimal convergence, that slowly deteriorates for larger penalty parameter, in agreement with the analysis in \cite{Brenner2005}. Variations of C0-IPM have also been applied to strain gradient dependent damage models in \cite{Wells2003} and to the Cahn-Hilliard equation in \cite{Engel2002,Brenner2012}. 

This paper develops the C0-IPM method for flexoelectricity and, as a simplified case, for strain gradient elasticity, for 2D and 3D computations. Section \ref{sec:problemStatement} presents the problem statement and recalls the weak form in $\mathcal{H}^2$. The derivation of the C0-IPM method for $\mathcal{C}^0$ approximations, not in $\mathcal{H}^2$, is presented in section \ref{sec:COIPM}. An eigenvalue problem to determine a large enough penalty parameter, ensuring coercivity of the strain gradient bilinear form, is derived in section \ref{sec:eigs}. Finally, in section \ref{sec:examples}, 2D and 3D numerical experiments  demonstrate the applicability of the method and show, as expected, slightly suboptimal convergence under uniform mesh refinement, but still with a robust high-order convergence for $p\geq 3$.

\section{Problem statement}
\label{sec:problemStatement}
We consider the model in \cite{Codony2019}, where flexoelectricity is ruled by the following set of PDEs and boundary conditions:
\begin{subequations} \label{eq:problemStatement}
\begin{align} 
\grad\cdot (\hsigma(\bu,\phi)-\grad\cdot \tsigma(\bu,\phi)) + \bm{b} = \bm{0} &\text{ in }\Omega \label{eq:mechanical}\\
    \nabla \cdot \hD(\bu,\phi) - q = 0&\text{ in }\Omega  \label{eq:electical}\\
         \bu = \bm{g}_1 &\text{ on } \overline{\Gamma^{\bu}_{D_1}} \label{eq:Dirichlet1}\\
\bm{t} (\bu,\phi)  = \bm{t}_n & \text{ on }  \Gamma^{\bu}_{N_1} \label{eq:Neumann1}\\
     \cfrac{\partial \bu}{\partial \bm{n}} = \bm{g}_2 &\text{ on } \Gamma^{\bu}_{D_2} \label{eq:Dirichlet2}\\
\bm{r} (\bu,\phi) = \bm{r}_n & \text{ on }  \Gamma^{\bu}_{N_2} \label{eq:Neumann2}\\
	 \bm{j}(\bu,\phi) = \bm{j}^{ext} & \text{  on } C_N^{\partial \Omega} \label{eq:cornersOmega} \\
    \phi = g_3 &\text{ on }\Gamma_D^\phi  \label{eq:Dirichlet3}\\
    w(\bm{u},\phi)=w_n &\text{ on }\Gamma_N^\phi \label{eq:Neumann3}
\end{align}
\end{subequations}
where $\Omega\subset \mathbb{R}^\nsd$ is the domain, the displacement $\bu$ and the electric potential  $\phi$ are the unknowns, $\hsigma$ and $\tsigma$ are the local and double stress tensors and $\hD$ is the electric displacement tensor, that is, 
\begin{align*}
    \hsigma = \mathbb{\bm{C}} : \bm{\varepsilon} - \bm{E} \cdot \bm{e} \quad  
   &\equiv \quad  \widehat{\sigma}_{ij} = C_{ijk\ell} \, \varepsilon_{k\ell} -  E_\ell\,e_{\ell ij} \\
      \tsigma = \bm{h} \, \vdots\, \bm{\nabla}\bm{ \varepsilon}-\bm{E}\cdot \bm{\mu}  \quad
     & \equiv \quad  \widetilde{\sigma}_{ijk} = h_{ijk\ell mn} \, \cfrac{\partial\varepsilon_{\ell m}}{\partial x_n} - E_\ell\,\mu_{\ell ijk}\\
         \widehat{\bm{D}} = \bm{\kappa} \cdot \bm{E}  +\bm{e}: \bm{\varepsilon} +\bm{\mu} \vdots \bm{\nabla} \bm{\varepsilon} \quad
         &\equiv \quad \widehat{D}_\ell = \kappa_{\ell m} E_m + e_{\ell ij}\varepsilon_{ij}+ \mu_{\ell ijk}\cfrac{\partial \varepsilon_{ij} }{\partial x_k}
\end{align*}
$\bm{\varepsilon}$ is the strain tensor, that is $\bm{\varepsilon}_{ij}=(\partial u_i / \partial x_j + \partial u_j / \partial x_i)/2$, $\bE=-\grad \phi$ is the electric field, 
 $\bm{C}$ is the elasticity tensor, that depends on the Young modulus $E$ and the Poisson ratio $\nu$, $\bm{h}$ is the strain-gradient tensor, defined as
$ h_{ijk\ell mn}=l^2 C_{ij\ell m}\delta_{kn}$ 
 with the internal length scale {parameter} $l$,
$\bm{e}$ and  $\bm{\mu}$ are the tensors of piezoelectric and flexoelectric coefficients  and $\bm{\kappa}$ contains the dielectricity constants, see appendix B in \cite{Codony2019} for detailed definitions.

In the previous equations, and in the rest of the document,
 Einstein's notation is assumed. That is, repeated indexes sum over the spatial dimensions. 
 
The boundary of the domain is split in Dirichlet and Neumann boundaries, for the first and second conditions of the mechanical problem and for the electric problem, that is
\[
\partial\Omega = \overline{\Gamma^{\bu}_{D_1} \cup  \Gamma^{\bu}_{N_1}} =  \overline{\Gamma^{\bu}_{D_2} \cup  \Gamma^{\bu}_{N_2}}  =  \overline{\Gamma^{\phi}_{D} \cup  \Gamma^{\phi}_{N}}.
\]
Note that all volume and boundary domains are assumed to be open domains, not including their boundaries.

The first mechanical boundary condition, \eqref{eq:Dirichlet1}  or \eqref{eq:Neumann1}, sets the displacement $\bu$ or the traction
\[
        t_i (\bm{u},\phi)= \left(\widehat{\sigma}_{ij} -\cfrac{\partial \widetilde{\sigma}_{ijk}}{\partial x_k} - \nabla^{S}_k \widetilde{\sigma}_{ikj}\right) n_j + \widetilde{\sigma}_{ijk } \widetilde{N}_{jk},
\]
where {$\nabla^{S}_k \widetilde{\sigma}_{ikj}$ is the surface divergence of $\widetilde{\sigma}_{ikj}$,} and $\tilde{N}$ is the second order geometry tensor, see \cite{Codony2019} for details. The second mechanical boundary condition, \eqref{eq:Dirichlet2}  or \eqref{eq:Neumann2}, sets the normal derivative of the displacement ${\partial \bu}/{\partial \bm{n}} $ or the double traction 
\[
    r_i(\bm{u},\phi)=\widetilde{\sigma}_{ijk}n_jn_k.
\]
\begin{figure}[t]
\centering
  \includegraphics[width=0.35\textwidth]{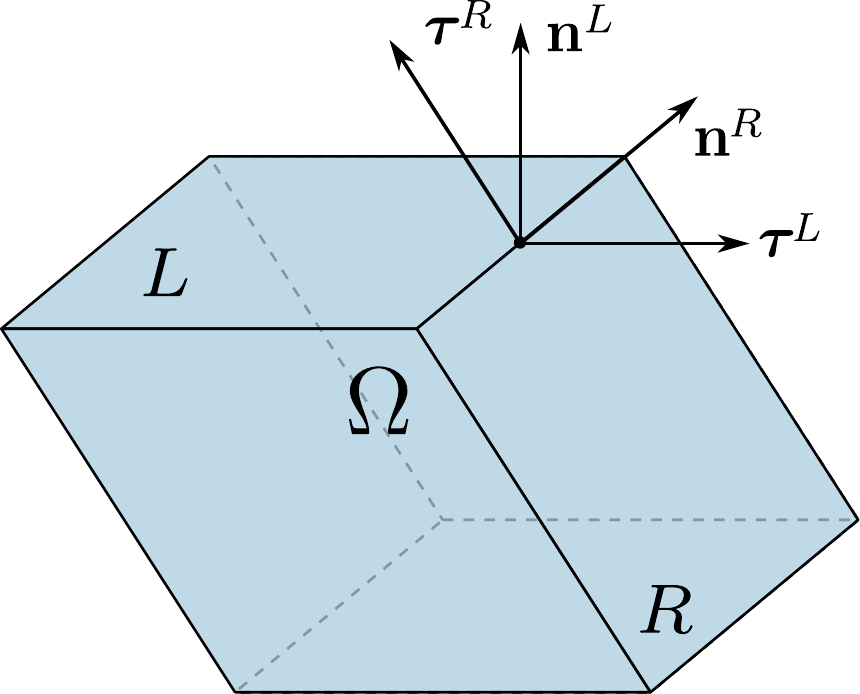}
  \hspace{1.5cm}
  \includegraphics[width=0.3\textwidth]{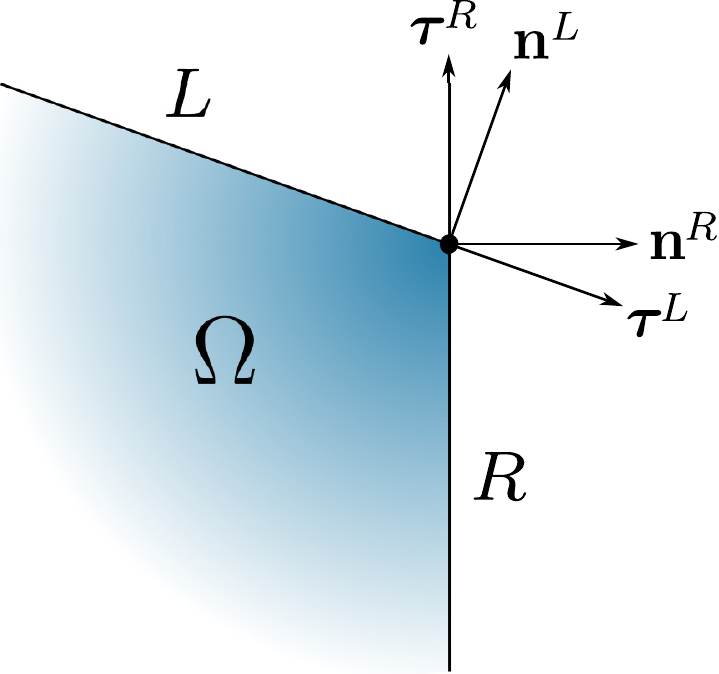}
\caption{{Sketch of normal and tangent vectors on an edge (3D, left) and at a corner (2D, right), for the computation of the corresponding line and punctual forces $\bm{j}(\bu,\phi)$. The superscripts $L$ and $R$ refer to the face (side in 2D) sharing the edge (corner in 2D).}}
\label{fig:cornerVectors}       
\end{figure}
The condition \eqref{eq:cornersOmega} sets forces on the Neumann boundary edges. That is, the domain boundary  is assumed to be composed of smooth surfaces (curves in 2D) that are joined on sharp boundary edges (corners in 2D). $C_N^{\partial \Omega}$ denotes the union of the boundary edges that are shared by two surfaces with first Neumann conditions, i.e. the edges in the interior of $\Gamma^{\bu}_{N_1}$. At the edges shared by at least one Dirichlet surface the value on the edge is assumed to be the one set on the surface, i.e. $\bu=\bm{g}_1$ for all edges in $\overline{\Gamma^{\bu}_{D_1}}$. Line forces (punctual forces in 2D) are defined on boundary edges as
\[
   j_i(\bm{u},\phi) = \tau_j^L \widetilde{\sigma}_{ij\ell}^Ln_\ell^L + \tau_j^R \widetilde{\sigma}_{ij\ell}^Rn_\ell^R,
\]
being $\bm{n}^L$ and $\bm{n}^R$  the unitary exterior normals on the left and right surfaces sharing the boundary edge, and $\bm{\tau}^L$ and $\bm{\tau}^R$ the tangent vectors on each surface pointing outward and perpendicular to the edge{, see an example in figure \ref{fig:cornerVectors} left. In 2D, $\bm{\tau}^L$ and $\bm{\tau}^R$ at a corner are just the tangent vectors on each curve sharing the corner and pointing outward, as depicted in figure \ref{fig:cornerVectors} right.}

Finally, the electric boundary condition,  \eqref{eq:Dirichlet3}  or \eqref{eq:Neumann3}, sets the electric potential $\phi$ or the surface charge density
\[
w(\bm{u},\phi) =-D_\ell(\bm{u},\phi)\, n_\ell.
\]

\begin{remark}
For the sake of simplicity, we initially restrict to the case with second Neumann boundary conditions in the whole boundary, i.e. $\Gamma^{\bu}_{N_2}=\partial \Omega$ and $\Gamma^{\bu}_{D_2}=\emptyset$. The treatment of second Dirichlet conditions \eqref{eq:Dirichlet2} is commented in Remark {\ref{rem:Dirichlet2}}.
It is also worth mentioning that the rationales in this work can also be applied to other models, including converse flexoelectricity or expressed in terms of polarization, see for instance \cite{Nanthakumar2017}.  
\end{remark}

If an approximation in $\mathcal{H}^2(\Omega)$ can be considered, multiplying  \eqref{eq:mechanical} by a weighting vector $\bv$, applying integration by parts twice, and using the symmetries of the stress tensors, leads to
\begin{align*}
   \int_{\Omega} \bv \cdot \bm{b} \, \text{d}\Omega = &\int_\Omega \bm{\varepsilon}(\bv) : {\hsigma}(\bu,\phi) \, \text{d}\Omega + \int_{\Omega} \grad \bm{\varepsilon}(\bv) \vdots \tsigma (\bu,\phi)\, \text{d}\Omega \\
   &- \int_{\partial \Omega }v_i \left(\widehat{\sigma }_{ij}-\cfrac{\partial  \widetilde{\sigma}_{ijk}}{\partial x_k}\right)\, n_j \, \text{dS} - \int_{\partial \Omega} \cfrac{\partial v_i}{\partial x_j}\, \widetilde{\sigma}_{ijk}\, n_k \, \text{dS},
 \end{align*}
 where $\grad \bm{\varepsilon} \vdots\tsigma=\partial \bm{\varepsilon}_{ij}/\partial x_k \,\widetilde{\sigma}_{ijk}$.
 
Now, to properly treat boundary conditions, the derivative ${\partial v_i}/{\partial x_j}$ on the boundary is split in normal and tangential derivatives, and the surface divergence theorem is applied to the term with the tangential derivative, leading to 
\begin{equation}\label{eq:weakFormH2}
\begin{split}
    &\int_{\Omega} \bv \cdot  \bm{b} \, \text{d}\Omega = \int_\Omega \bm{\varepsilon}(\bv) : {\hsigma}(\bu,\phi) \, \text{d}\Omega + \int_{\Omega}\grad \bm{\varepsilon}(\bv) \vdots \tsigma (\bu,\phi) \, \text{d}\Omega  \\
    &- \int_{\partial \Omega} 	\bv \cdot \bm{t}( \bm{u},\phi) \, \text{dS}- \int_{\partial\Omega}  \cfrac{\partial \bv }{\partial\bm{n}} \cdot \bm{r}(\bm{u},\phi)\, \text{dS} - \int_{C^{\partial\Omega}}\bv \cdot \bm{j}(\bm{u},\phi)  \,\text{dl},
    \end{split}
\end{equation}
where $C^{\partial\Omega}$ is the union of all sharp edges of the domain, and the integral on it reduces to a sum evaluating at the boundary corners in 2D.

Thus, applying boundary conditions \eqref{eq:Dirichlet1}-\eqref{eq:Neumann3}, under the assumption $\Gamma^{\bu}_{D_2}=\emptyset${, and adding the weighted residual of the electric potential problem \eqref{eq:electical} with \eqref{eq:Dirichlet3} and \eqref{eq:Neumann3}}, the weak form of  \eqref{eq:problemStatement} in  $\mathcal{H}^2(\Omega)$ is: find $\bu\in [\mathcal{H}^2(\Omega)]^\nsd$ and $\phi \in \mathcal{H}^1(\Omega)$ such that \eqref{eq:Dirichlet1} and  \eqref{eq:Dirichlet3} hold and
\begin{equation}
\label{eq:wfH2}
\begin{split}    
    \int_\Omega \bm{\varepsilon}(\bv) : {\hsigma}(\bu,\phi) \, \text{d}\Omega &+ \int_{\Omega} \grad\bm{\varepsilon}(\bv) \vdots \tsigma (\bu,\phi) \, \text{d}\Omega  + \int_\Omega \grad \psi \cdot \hD(\bu,\phi)\, \text{d}\Omega  \\
    &= s(\bv)
    - \int_{\Omega} \psi q \;\text{d}\Omega -\int_{\Gamma^{\phi}_N}\psi\, w_n\text{ d}S
\end{split}
\end{equation}
for all $\bv\in [\mathcal{H}^2(\Omega)]^\nsd$ and $\psi \in \mathcal{H}^1(\Omega)$ such that $\bv=\bm{0}$ on $\overline{\Gamma^{\bu}_{D_1}}$ and $\psi=0$ on $\Gamma^{\phi}_{D}$ where
\begin{equation} \label{eq:sdef}
s(\bv)= \int_{\Gamma^{\bu}_{N_1}} \bv \cdot \bm{t}_n \, \text{dS}+\int_{\Gamma^{\bu}_{N_2} } \cfrac{\partial \bv }{\partial\bm{n}}\cdot \bm{r}_n\, \text{dS} + \int_{C_N^{\partial\Omega}}\bv \cdot \bm{j}^{ext} \,\text{dl}  + \int_{\Omega} \bv \cdot \bm{b} \, \text{d}\Omega.
\end{equation}

This weak form is not suitable when considering $\mathcal{C}^0$ FE approximations, but the same derivation can be applied in the interior of each element of the mesh, as detailed next.

\section{$\mathcal{C}^0$ Interior Penalty Finite Element method}
\label{sec:COIPM}

The domain $\Omega$ is now split in FEs $\{\Omega_e\}_{e=1}^\numel$, and a $\mathcal{C}^0$ piece-wise polynomial approximation is considered. That is, the approximation space {for the components of the displacement and for the potential} is
\[
\mathcal{V}^h=\{v\in \mathcal{H}^1(\Omega) \text{ such that } \varphi_e^{-1}(v|_{\Omega_e} )\in \mathcal{P}^p \text{ for } e=1,\dots,\numel\} \not\subset \mathcal{H}^2(\Omega),
\]
where $\varphi_e$ is the isoparametric transformation from the reference element to the physical element $\Omega_e$, and $\mathcal{P}^p$ is the space of polynomials of degree less or equal to $p$ for simplexes, and less or equal to  $p$  in each direction for quadrilaterals and hexahedra.  

\begin{figure}[!htb]
\centering
  \includegraphics[width=\textwidth]{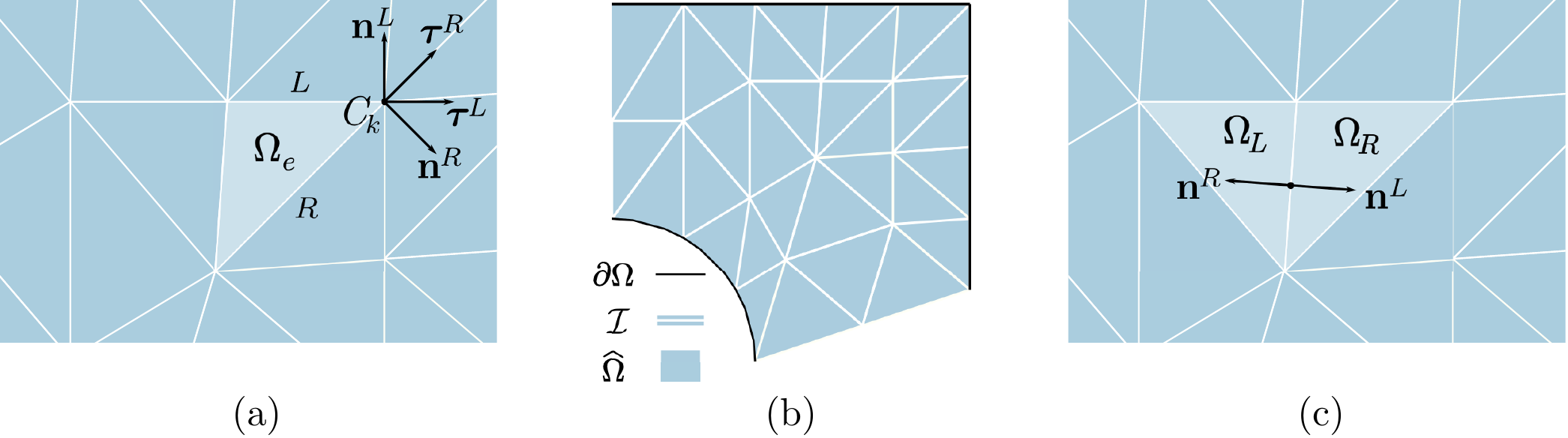}
\caption{Example of 2D discretization: (a) interior corner $C_k$ (i.e. interior mesh vertex), as seen from element $\Omega_e$, and representation of the normal and tangent vectors corresponding to the left and right side of the element sharing the corner, for the definition of the punctual force $\bm{j}^e$ at the corner, (b) interior faces $\I$ in white and broken domain $\widehat\Omega$ in blue, (c) normal vectors on one face shared by its left and right elements, for the computation of the jump on a face.}
\label{fig:tangentsNormals}       
\end{figure}

Since the approximation space in not in $\mathcal{H}^2(\Omega)$, we can not consider the weak form \eqref{eq:wfH2}. However, the approximation is $\mathcal{H}^2(\Omega_e)$; thus, considering 
\eqref{eq:weakFormH2} in each element we have
\begin{equation}\label{eq:elements}
\begin{split}
  &  \int_{\Omega_e} \bv \cdot  \bm{b} \, \text{d}\Omega = \int_{\Omega_e} \bm{\varepsilon}(\bv) : {\hsigma}(\bu,\phi) \, \text{d}\Omega + \int_{\Omega_e} \grad\bm{\varepsilon}(\bv) \vdots \tsigma (\bu,\phi) \, \text{d}\Omega 
    \\  &
    - \int_{\partial \Omega_e} \bv \cdot \bm{t}^e( \bm{u},\phi) \, \text{dS} - \int_{\partial\Omega_e}  \cfrac{\partial \bv }{\partial\bm{n}^e} \cdot \bm{r}^e(\bm{u},\phi)\, \text{dS} - \int_{C^{\partial\Omega_e}}\bv \cdot \bm{j}^e(\bm{u},\phi) \,\text{dl},
\end{split}
\end{equation}
where $C^{\partial\Omega_e}$ is the union of the edges (corners in 2D) of the element $\Omega_e$ and $\bm{j}^e$ is the line force (punctual force in 2D) on $C^{\partial\Omega_e}$; see a representation in figure \ref{fig:tangentsNormals}a. The superscripts highlight that  the surface and line forces, and the normal vector, are from the element $\Omega_e$.

Summing  for all elements, and noting that $\bv$ is continuous but  ${\partial \bv }/{\partial\bm{n}}$ is not, we get
\begin{equation}\label{eq:elementsSum}
\begin{split}
    \int_{\Omega} \bv \cdot \bm{b} \, \text{d}\Omega = &\int_{\Omega} \bm{\varepsilon}(\bv) : {\hsigma}(\bu,\phi) \, \text{d}\Omega + \int_{\widehat\Omega} \grad \bm{\varepsilon}(\bv) \vdots \tsigma (\bu,\phi) \, \text{d}\Omega 
    \\
    &- \int_{\I} \bv \cdot(\bm{t}^L( \bm{u},\phi)+\bm{t}^R( \bm{u},\phi)) \, \text{dS} 
    \\
    &-\int_{\I}   \left( \cfrac{\partial \bv^L }{\partial\bm{n}^L}\cdot \bm{r}^L(\bm{u},\phi)+ \cfrac{\partial \bv^R }{\partial\bm{n}^R}\cdot \bm{r}^R(\bm{u},\phi)\right) \text{dS} 
    \\
    &-\sum_{k=1}^{\nedges} \int_{C_k} \bv \cdot \left(\sum_{e\in E(k)} \bm{j}^e(\bm{u},\phi) \right) \text{dl}
    \\
 &-\int_{\partial \Omega} \bv \cdot \bm{t}( \bm{u},\phi) \, \text{dS} 
    - \int_{\partial\Omega}  \cfrac{\partial \bv }{\partial\bm{n}} \cdot \bm{r}(\bm{u},\phi)\, \text{dS},   
\end{split}
\end{equation}
where $\I$ is the union of all the internal element faces and $\widehat\Omega$ is the union of the interior of the elements, {where second derivatives are well-defined,} i.e. 
\[
\widehat\Omega= \bigcup_{e=1}^\numel \Omega_e, \quad \I=\bigcup_{e=1}^\numel \partial \Omega_e \backslash \partial \Omega,
\]
see figure \ref{fig:tangentsNormals}b.
The supercripts $R$ and $L$ now denote the evaluation from the elements to the left and right side of the face in $\I$ (see figure \ref{fig:tangentsNormals}c), and 
$\{C_k\}_{k=1}^{\nedges}$ are all the  edges (corners in 2D) in the mesh, being  $E(k)$ the set of indexes of the elements sharing the edge $C_k$.

Now, let us recall that 
the conditions for interfaces in the domain (also in the case of discontinuous material parameters) are the ones corresponding to both continuity for the Dirichlet values and equilibrium of Neumann forces. That is, 
\begin{subequations} \label{eq:interfaceConditions}
\begin{align}
         &\jump{\bu \otimes \bm{n}} = \bm{0}  \label{eq:Dirichlet1_I}\\
&\jump{\bm{t} (\bu,\phi)}  = \bm{0}   \label{eq:Neumann1_I}\\
     &\jump{\cfrac{\partial \bu}{\partial \bm{n}} }= \bm{0} \label{eq:Dirichlet2_I}\\
&\jump{\bm{r} (\bu,\phi)\otimes \bm{n}}  =  \bm{0} \label{eq:Neumann2_I}
\end{align}
\end{subequations}
on the faces in $\I$, where
the jump operator is defined as
\begin{align*}
  \jump{a} = a^L+a^R
\end{align*}
and it is used always involving a change of sign due to an odd appearance of the normal vector.

In addition, we have to impose equilibrium of forces on the mesh edges.
That is, each element $\Omega_e$ contributes with a force $\bm{j}^e$ on its edges (corners in 2D); and for each edge $C_k$, the sum of the forces for all elements sharing the edge (i.e. for $\Omega_e$ with $e\in E(k)$) must be zero, or in internal equilibrium with the external forces. That is,
\begin{equation}
\sum_{e\in E(k)}  \bm{j}^e(\bm{u},\phi) =  
\left\{\begin{array}{ll}
\bm{0} &\text{ on } C_k \not\subset C^{\partial \Omega} \\
\bm{j}^{ext} &\text{ on } C_k \subset C_N^{\partial \Omega}
\end{array}
\right.
\label{eq:cornersI} 
\end{equation} 
where $\bm{j}^{ext}$ is the force set in  \eqref{eq:cornersOmega}.
Note that $\{C_k \not\subset C^{\partial \Omega}\}$ includes interior edges and also element edges on $\partial \Omega$, just excluding the ones in the domain sharp edges.
For the edges in $\overline{\Gamma_{D_1}^{\bu}}$ no value is set, and the sum of the forces will be in equilibrium with the reaction forces associated to the prescribed displacement \eqref{eq:Dirichlet1}.

On other hand, using the algebraic identity $(a^L\bm{n}^L)b^L+(a^R\bm{n}^R)b^R = \mean{a} \jump{b\bm{n}} + \jump{a\bm{n}}\mean{b}$, and the equilibrium condition \eqref{eq:Neumann2_I}, we can rewrite 
\begin{equation}\label{eq:vnr}
\begin{split}
    \cfrac{\partial \bv^L }{\partial\bm{n}^L}\cdot \bm{r}^L(\bm{u},\phi)&+ \cfrac{\partial \bv^R }{\partial\bm{n}^R}\cdot \bm{r}^R(\bm{u},\phi) 
    \\ &
    = \mean{\bm{\nabla v}} : \jump{\bm{r}(\bm{u},\phi)\otimes \bm{n}} + \jump{\cfrac{\partial \bm{v}}{\partial \bm{n}}} \cdot\mean{\bm{r}(\bm{u},\phi)} 
    \\ & 
    = \jump{\cfrac{\partial \bm{v}}{\partial \bm{n}}} \cdot\mean{\bm{r}(\bm{u},\phi)},
\end{split}
\end{equation}
with the mean operator $\{a\} = (a^L+a^R)/2$.

Now, replacing in \eqref{eq:elementsSum} the identity \eqref{eq:vnr}, the Neumann boundary conditions \eqref{eq:Neumann1} and \eqref{eq:Neumann2}, the homogeneous Dirichlet condition $\bv=\bm{0}$ on $\overline{\Gamma_{D_1}^{\bu}}$ related to \eqref{eq:Dirichlet1}, the first interface equilibrium condition \eqref{eq:Neumann1_I} and the equilibrium at interior edges \eqref{eq:cornersI}, and under the assumption $\Gamma^{\bu}_{D_2}=\emptyset$, \eqref{eq:elementsSum} simplifies to
\begin{equation} \label{eq:wfH1_1}
\begin{split}
   \int_{\Omega} \bm{\varepsilon}(\bv) : {\hsigma} (\bu,\phi) \, \text{d}\Omega &+ \int_{\widehat\Omega} \grad\bm{\varepsilon}(\bv) \vdots \tsigma (\bu,\phi) \, \text{d}\Omega 
   \\   &
   -\int_{\I}   \jump{\cfrac{\partial \bm{v}}{\partial \bm{n}}}\cdot \mean{\bm{r}(\bm{u},\phi)} \text{dS}
    = s(\bv),  
\end{split}
\end{equation}
with $s(\bv)$ defined in \eqref{eq:sdef}.

The first two integrals in \eqref{eq:wfH1_1} are symmetric and coercive bilinear forms in $\bv$ and $\bu$, as expected for the weak form of a strain gradient elasticity operator. However, it is not the case for the integral on the interior faces $\I$.

The idea of IPM is adding terms that are analytically zero, thanks to the continuity interface condition \eqref{eq:Dirichlet2_I}, to recover symmetry and coercivity of the strain gradient bilinear form. 
The resulting weak form for flexoelectricity, under the assumption $\Gamma^{\bu}_{D_2}=\emptyset$, is: find $\bu\in [\mathcal{H}^1(\Omega)\cap \mathcal{H}^2(\widehat\Omega)]^\nsd$ and $\phi \in \mathcal{H}^1(\Omega)$ such that \eqref{eq:Dirichlet1} and  \eqref{eq:Dirichlet3} hold and
\begin{equation}
\label{eq:wfH1}
\begin{split}    
    \int_\Omega \bm{\varepsilon}(\bv) : {\hsigma}(\bu,\phi) & \text{d}\Omega +\int_{\widehat\Omega} \grad\bm{\varepsilon}(\bv) \vdots \tsigma (\bu,\phi) \, \text{d}\Omega 
    + \int_\Omega \grad \psi \cdot \hD(\bu,\phi) \, \text{d}\Omega 
    \\
    &-\int_{\I}   \jump{\cfrac{\partial \bm{v}}{\partial \bm{n}}}\cdot \mean{\bm{r}(\bm{u},\phi)} \text{dS}
     -\int_{\I}   \mean{\bm{r}(\bm{v},\psi)} \cdot \jump{\cfrac{\partial \bm{u}}{\partial \bm{n}}} \text{dS}
     \\
    &+\int_{\I}  \beta \jump{\cfrac{\partial \bm{v}}{\partial \bm{n}}} \cdot \jump{\cfrac{\partial \bm{u}}{\partial \bm{n}}} \text{dS}
    \\
    &= s(\bv)- \int_{\Omega} \psi q \;\text{d}\Omega -\int_{\Gamma^{\phi}_N}\psi\, w_n\text{ d}S
\end{split}
\end{equation}
for all $\bv\in [\mathcal{H}^1(\Omega)\cap \mathcal{H}^2(\widehat\Omega)]^\nsd$ and $\psi \in \mathcal{H}^1(\Omega)$ such that $\bv=\bm{0}$ on $\Gamma^{\bu}_{D_1}$ and $\psi=0$ on $\Gamma^{\phi}_{D}$.

The parameter $\beta$ is a stabilization parameter that must be taken large enough to ensure coercivity of the strain gradient bilinear form, to get a well-defined saddle point problem \cite{Codony2019}. Although it is usually called penalty parameter, thanks to the consistency of the formulation, high-order convergence can be achieved with $\beta$ or order $h^{-1}$. In practice, not very large values are needed for accurate solutions, avoiding the unaccuracy or ill-conditioning that typically suffer non-consistent penalty methods \cite{Fernandez2004}.

The minimum value of the stabilization parameter $\beta$ can be estimated solving an eigenvalue problem, as commented in section \ref{sec:eigs}.

\begin{remark}[Second Dirichlet conditions] \label{rem:Dirichlet2}
If second Dirichlet boundary conditions \eqref{eq:Dirichlet2} are imposed (i.e. ${\Gamma_{D_2}^{\bu}}\neq \emptyset$), an additional term $\int_{\Gamma_{D_2}^{\bu}}{\partial \bm{v}}/{\partial \bm{n}}\cdot \bm{r}(\bm{u},\phi) \,\text{dS}$ appears in \eqref{eq:wfH1_1} and, consequently, in the weak form \eqref{eq:wfH1}. Following the same IPM rationale, two new terms, that are null thanks to \eqref{eq:Dirichlet2_I}, are also added in \eqref{eq:wfH1} to recover again symmetry and coercivity, namely 
$\int_{\Gamma_{D_2}^{\bu}} \bm{r}(\bm{v},\psi) \cdot \left({\partial \bm{u}}/{\partial \bm{n}}- \bm{g}_2 \right) \, \text{dS} 
+
\beta_D \int_{\Gamma_{D_2}^{\bu}} {\partial \bm{v}}/{\partial \bm{n}} \cdot \left({\partial \bm{u}}/{\partial \bm{n}}- \bm{g}_2 \right) \, \text{dS} $, where $\beta_D$ is a new stabilization parameter that can be taken equal to $\beta$ or tuned separately.
\end{remark}

\begin{remark}
It is interesting to note that the C0-IPM weak form \eqref{eq:wfH1} reduces to the one for $\mathcal{H}^2(\Omega)$, i.e.  \eqref{eq:wfH2}, when a  $\mathcal{C}^1(\Omega)$ approximation is considered. The C0-IPM formulation keeps the consistency and is valid for standard $\mathcal{C}^0$ FE approximations, just introducing the proper integrals on the faces $\I$. Also note that the second integral in \eqref{eq:wfH1} is in the interior of the elements, $\widehat \Omega$, to account for the fact that second derivatives are not defined on $\I$. 
\end{remark}

\subsection{Estimate of the interior penalty parameter $\beta$}
\label{sec:eigs}

In this section we derive an eigenvalue problem to {estimate a lower bound for} $\beta$. The derivation is the usual one in IPM and Nitsche's formulations \cite{Griebel2002,Codony2019,Fojo2020}.

The bilinear form of the strain-gradient elasticity operator is 
\begin{equation*}
\begin{split}
\mathcal{A}(\bu,\bv) = a(\bu,\bv) 
    &-\int_{\I}   \jump{\cfrac{\partial \bm{v}}{\partial \bm{n}}}\cdot \mean{\bm{r}^{sg}(\bm{u})} \text{dS}
     -\int_{\I}   \mean{\bm{r}^{sg}(\bm{v})} \cdot \jump{\cfrac{\partial \bm{u}}{\partial \bm{n}}} \text{dS}
     \\&
    + \beta \int_{\I}   \jump{\cfrac{\partial \bm{v}}{\partial \bm{n}}} \cdot \jump{\cfrac{\partial \bm{u}}{\partial \bm{n}}} \text{dS}
\end{split}
\end{equation*}
with 
\[
a(\bu,\bv) = 
\int_{\widehat\Omega} \grad\bm{\varepsilon}(\bv) \vdots \tsigma^{sg} (\bu) \, \text{d}\Omega 
\]
where $r^{sg}_i=\widetilde{\sigma}^{sg}_{ijk}n_jn_k$
and $\tsigma^{sg}(\bu)=\bm{h} \, \vdots\, \bm{\nabla}\bm{ \varepsilon}(\bu)$, that is, the mechanical part of the second traction and the double stress tensor.

The bilinear form $a$ is semicoercive (i.e. $a(\bu,\bu)>0$ for any $\bu$ such that $\grad \bm\varepsilon (\bu) \neq \bm{0}$, and $a(\bu,\bu)=0$ otherwise), leading to {well-posed strain gradient elasticity  and flexoelectricity  problems} for any value of the {internal length scale parameter} $l$. However, the addition of the integrals on the faces $\I$, leads to a bilinear form $\mathcal{A}$ that retains semicoercivity only for large enough $\beta$. 

Thus, to ensure well-posedness of the discrete problem for any value of $l$, we want $\beta$ such that
$\mathcal{A}(\bu,\bu)>0  \; \forall \; \bu\in \widetilde{\mathcal{U}}^h=\{\bu \in [\mathcal{V}^h]^\nsd \text{ such that } \grad \bm\varepsilon (\bu) \neq \bm{0} \}$ with 
\begin{equation*}
\mathcal{A}(\bu,\bu) = a(\bu,\bu) 
    -2\int_{\I}   \jump{\cfrac{\partial \bm{u}}{\partial \bm{n}}}\cdot \mean{\bm{r}^{sg}(\bm{u})} \text{dS}
    +\beta \int_{\I}   \jump{\cfrac{\partial \bm{u}}{\partial \bm{n}}} \cdot \jump{\cfrac{\partial \bm{u}}{\partial \bm{n}}} \text{dS}.
\end{equation*}

Using the Cauchy-Schwarz and Young's inequalities we can bound the interface terms as
\begin{equation*}
\begin{split}
    -2\int_{\I}   &\jump{\cfrac{\partial \bm{u}}{\partial \bm{n}}}\cdot \mean{\bm{r}^{sg}(\bm{u})} \text{dS}
    +\beta \int_{\I}   \jump{\cfrac{\partial \bm{u}}{\partial \bm{n}}} \cdot \jump{\cfrac{\partial \bm{u}}{\partial \bm{n}}} \text{dS}
\\&   
   \geq - 2\norm{\mean{\bm{r}^{sg}(\bm{u})} } \norm{\jump{\cfrac{\partial \bm{u}}{\partial \bm{n}}}} 
   + \beta \norm{\jump{\cfrac{\partial \bm{u}}{\partial \bm{n}}}}^2
 \\&   
   \geq   - \frac{1}{\epsilon} \norm{\mean{\bm{r}^{sg}(\bm{u})} } ^2
   + (\beta- \epsilon) \norm{\jump{\cfrac{\partial \bm{u}}{\partial \bm{n}}}}^2,
\end{split}
\end{equation*}
for any positive $\epsilon$.

Thus, considering a positive constant $K$ such that
\begin{equation}\label{eq:vap1}
\norm{\mean{\bm{r}^{sg}(\bm{u})} } ^2 \leq K a(\bu,\bu)
\quad \forall \bu \in \widetilde{\mathcal{U}}^h, 
\end{equation}
we have
\[
\mathcal{A}(\bu,\bu) \geq \left(1-\frac{K}{\epsilon}\right) a(\bu,\bu)+ (\beta- \epsilon) \norm{\jump{\cfrac{\partial \bm{u}}{\partial \bm{n}}}}^2
\]
and the bilinear form is then positive definite if $1-K/{\epsilon}>0$ and $\beta- \epsilon>0$, for any positive $\epsilon$.

In conclusion, the strain gradient  bilinear form $\mathcal{A}$ is positive definite in the reduced discrete space if $\beta>K$, where $K$ is the constant satisfying \eqref{eq:vap1}. 
This constant can be computed as the largest eigenvalue of the {generalised problem}
\[
\bf{B} \bf{x} = \lambda \bf{A} \bf{x}
\]
where $\bf{B}$ and $\bf{A}$ are the discrete matrices corresponding to 
the bilinear forms $b(\bu,\bv)$ and $a(\bu,\bv)$ in the reduced discrete space $\widetilde{\mathcal{U}}^h$, with
\[
b(\bu,\bv)= \int_{\I}   \mean{\bm{r}^{sg}(\bm{v})} \cdot \mean{\bm{r}^{sg}(\bm{u})} \text{dS}.
\] 

The computation of the maximum eigenvalue in the reduced space $\widetilde{\mathcal{U}}^h$ can be done from the problem stated in the complete discrete space $[\mathcal{V}^h]^\nsd$ setting nodal values to reduce the space or using the so-called eigenvalue problem deflation \cite{Annavarapu2012}.

\begin{remark}\label{rem:beta}
Matrices $\bm{B}$ and $\bm{A}$ scale as $\mathcal{O}(E^2 l^4 h^{\nsd-5})$ and $\mathcal{O}(E l^2 h^{\nsd-4})$, respectively, with characteristic element size $h$. Thus, the maximum eigenvalue of \eqref{eq:vap1} scales as $\mathcal{O}(E l^2/h)$. Consequently, we can consider
\begin{equation}\label{eq:beta}
\beta = \alpha El^2/h,
\end{equation}
with a large enough constant $\alpha$, that can be computed solving the eigenvalue problem, or simply tuned, in a coarse mesh with any value of $E$ and $l$. 

It is important noting that, differently to non-consistent penalty methods, IPM and Nitsche's methods provide accurate solutions and high-order convergence  with moderate values of $\beta$ of order  $\mathcal{O}(h^{-1})$ for any degree of approximation, as shown in the numerical examples in section \ref{sec:examples} and in \cite{Fernandez2004}. 
\end{remark}

An alternative sufficient condition to have a well-posed discrete problem can be stated  including also the elasticity term in the bilinear form, that is, with 
\[
a(\bu,\bv) = 
\int_\Omega \bm{\varepsilon}(\bv) : {\hsigma}^{sg}(\bu) \text{d}\Omega +
\int_{\widehat\Omega} \grad\bm{\varepsilon}(\bv) \vdots \tsigma^{sg} (\bu) \, \text{d}\Omega ,
\]
{where ${\hsigma}^{sg}(\bu)= \mathbb{\bm{C}} : \bm{\varepsilon}(\bu)$  is the mechanical part of the local stress tensor.}
This option leads to a smaller (sharper) bound for $\beta$, specially for small $l$ or large $h$. However, since the matrix corresponding to the first elasticity term scales as $\mathcal{O}(Eh^{\nsd-2})$,  the dependency on the mesh size and material parameters is not so obvious.
\subsection{Implementation aspects}
\label{sec:implementationAspects}

The current implementation considers high-order Lagrange nodal basis, with Fekete nodes in the reference element to minimise the condition number of elemental matrices. For degree $p\geq$ 3, special attention must be paid to the position of interior nodes in curved physical elements to keep high-order convergence, see \cite{Chen1995}. High-order mesh generators, see for instance \cite{RuizGirones2019}, produce curved meshes taking care of this important aspect.

The computation of the system involves two separated loops: in elements for volume integrals, and in faces for the computation of integrals on $\I$.
To do so, the standard $\mathcal{C}^0$ reference element is extended including second derivatives of the basis functions at the element integration points, the value of element basis functions and their derivatives at the integration points of the reference element faces, a list of the nodes corresponding to each face in the reference element and permutations for the integration points of the reference face for \emph{flipping}. 

The so-called \emph{flipping} is a permutation (usually for the nodes in DG methods, but for integration points in our IPM implementation) that has to be applied to {the face when seen from the second element,} to match the orientation  of the corresponding face in the first element. In 2D the flipping is the same for any side of the mesh, just using a reverse ordering for the second element sharing the side. In 3D the possible rotations of the face have to be taken into account to choose the proper permutation for the integration points. 

A variable storing, for each face, the number of the elements sharing the face, the local numbering of the face in each one of the two elements and the rotation to be applied for the second element, is also computed from the mesh as a preprocess.

Dirichlet conditions \eqref{eq:Dirichlet1} are imposed in strong form, just setting the corresponding nodal values, as usual in standard FE computations. Second Dirichlet conditions can be imposed in weak form as commented in Remark \ref{rem:Dirichlet2}.

\section{Numerical examples}
\label{sec:examples}

Several numerical examples are included in this section to study the convergence of the C0-IPM formulation in 2D and 3D, and to validate the computational tool by comparison with previous works. 
Homogeneous first, second and corner Neumann boundary conditions are assumed where no boundary condition is specified.

\subsection{2D convergence test}
\label{sec:2Dconvergence}

\begin{figure}[!b]
    \centering
    \includegraphics[width=0.4\textwidth]{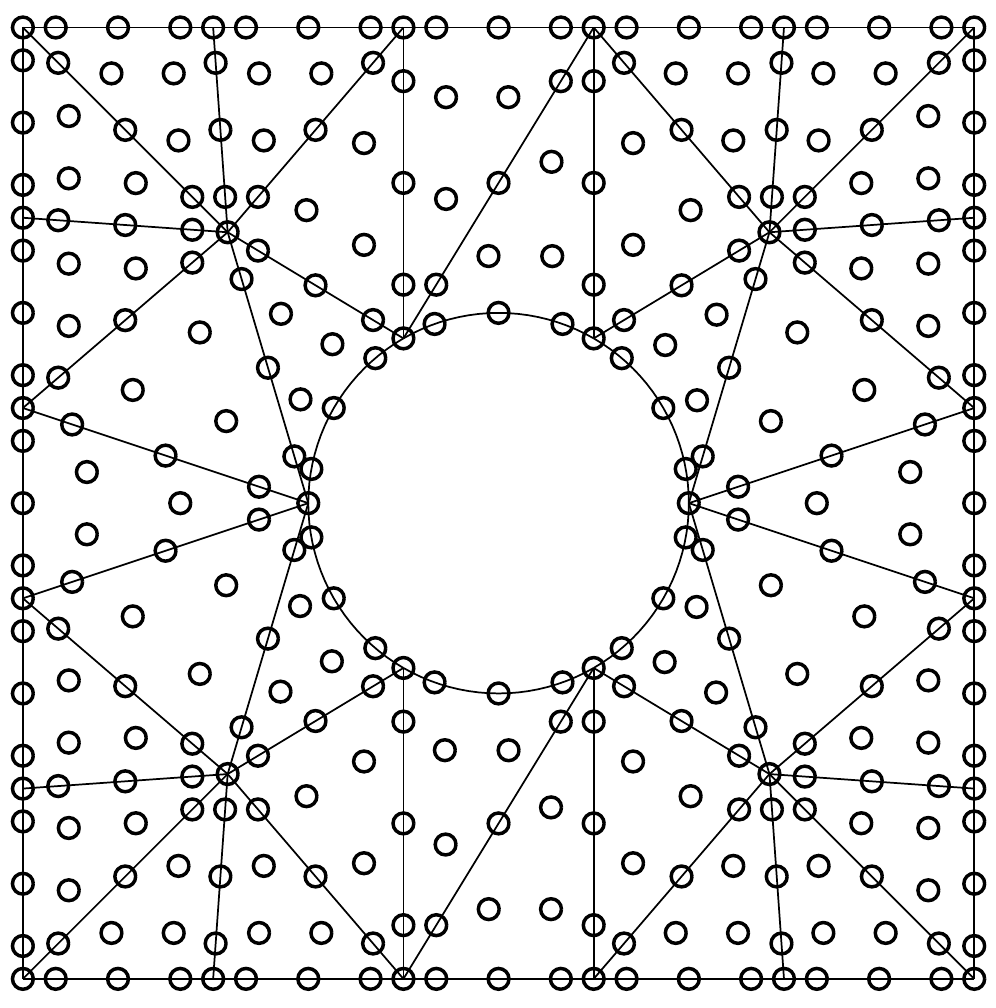}
    \caption{2D convergence test: {initial mesh for the nested refinement, with degree $p=4$.} }
    \label{fig:holeMesh}
\end{figure}

The convergence of the method for the solution of 
problem \eqref{eq:problemStatement} is studied in this section. To test the method with non-regular meshes and curved boundaries, the problem is solved in a square with a hole, $\Omega=(0,1)^2\backslash B((0.5,0.5),0.2)$. Figure \ref{fig:holeMesh} shows {the coarsest mesh for nested refinement, with degree $p=4$}.

\begin{figure}[!t]
    \centering
    \includegraphics[width=0.9\textwidth]{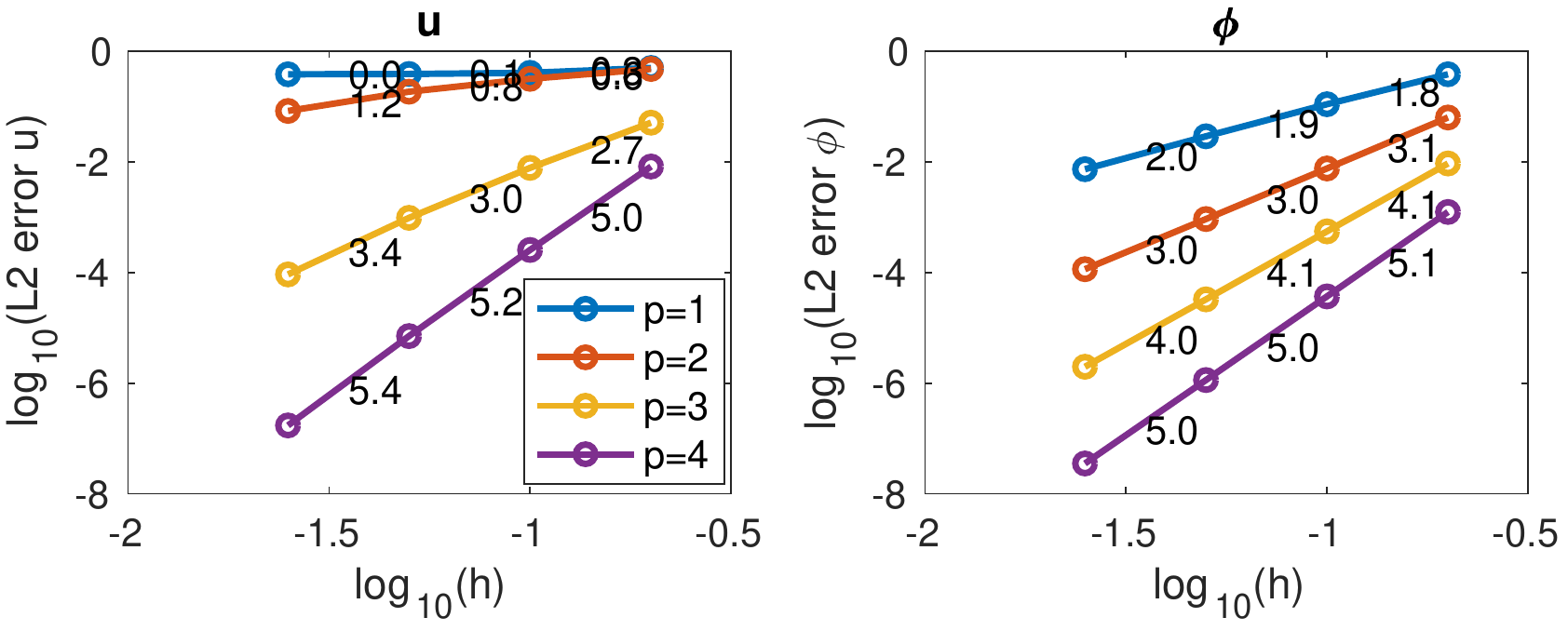}
    \caption{2D convergence test for the {uncoupled problem} (\emph{strain gradient elasticity} and potential equations) with  $\beta=100El^2/h$: $\mathcal{L}_2$ error under nested mesh refinement (the {coarsest mesh} is shown in figure \ref{fig:holeMesh}), for degree of approximation $p=1\dots 4$. The numbers are the slopes for each segment.}
    \label{fig:errors_Hole_alpha100_uncoupled}
    \includegraphics[width=0.9\textwidth]{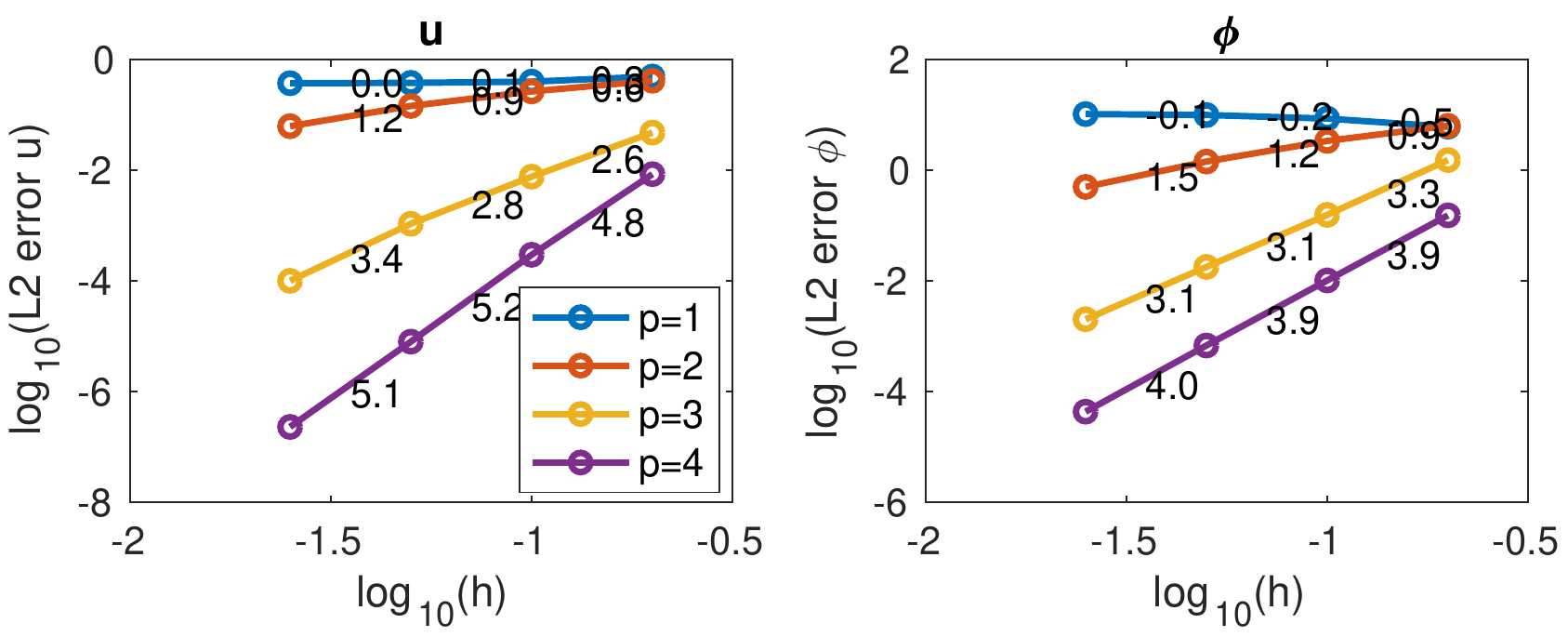}
    \caption{2D convergence test for the coupled \emph{flexoelectricity problem} with  $\beta=100El^2/h$. 
}
    \label{fig:errors_Hole_alpha100}
\end{figure}

First Dirichlet  and second Neumann conditions are imposed on all the boundary. The body force $\bm{b}$, the free charge $q$ and the boundary data are set so that the solution is 
\begin{equation}\label{eq:sinAnalytical}
\begin{array}{c}
	\bm{u} = \left[\sin{(2\pi(x_1+x_2))}, \cos{(2\pi(x_1+x_2))}\right]^T, \\
	  \phi = \sin(2\pi(x_1+x_2))+ \cos(2\pi(x_1+x_2)),
	 \end{array}
\end{equation}
the material parameters are
\begin{equation} \label{eq:materialParameters}
\begin{array}{lll}
     E = 2.5,      & \nu =0.25, & \\
     l=1.1,       & \kappa_L = 1.21, &\\
     e_L = 7.2
     , & e_T =1.33
     , & e_S = 1.73, \\
	  \mu_L = 1.5, & \mu_T = 1.34, & \mu_S = 5.47 ,
\end{array}
\end{equation}
and the piezoelectric principal direction is $x_1$. The definition of the material tensors in terms of these parameters can be found, for instance, in appendix B of \cite{Codony2019}.

First we consider the uncoupled problem (with $\bm{e}=\bm{0}$ and $\bm{\mu}=\bm{0}$), that is, an uncoupled solution of a strain gradient elasticity problem and an electric potential problem. The convergence plots are shown in figure \ref{fig:errors_Hole_alpha100_uncoupled}, for penalty parameter $\beta=100El^2/h$, and degree $p=1\dots 4$. For strain gradient elasticity, the displacement error behaves in agreement with the results for Kirchhoff plates in \cite{Fojo2020}. 
With degree $p=1$, the approximation space is not rich enough to impose $\mathcal{C}^1$ continuity. {Moreover, the second derivatives of the displacement in the strain gradient elasticity terms and the flexoelectricity terms cancel out, or are almost zero for curved elements. Thus, the method does not converge for linear approximation.} For degree $p= 2$ much finer meshes would be necessary to reach assymptotic convergence, reducing its practical applicability. Accurate results with high-order convergence are obtained for degree $p\geq 3$, with slightly suboptimal convergence for $p=3$, in agreement with the analysis in \cite{Brenner2005} for the biharmonic equation. In this particular example, $p=4$ behaves better than expected, exhibiting slightly superoptimal convergence. The expected optimal convergence is observed for the uncoupled electric potential problem for any degree.  

Figure \ref{fig:errors_Hole_alpha100} shows the convergence plots for the flexoelectricity problem, with piezoelectric and flexoelectric coupling. The coupling leads to a {reduction} in {the convergence rate}, not relevant for the displacement, but around {one} for the potential, for $p\geq 3$. This is probably due to the, small but still present, discontinuity of the displacement derivative across element sides, affecting the potential through the flexoelectricity coupling.

The conclusion is then that, even though convergence is suboptimal, the method is able to reach high accuracy with high-order convergence for degree $p\geq 3$. The C0-IPM method is therefore promising for an efficient solution of flexoelectricity.

Similar results can be observed with quadrilateral meshes, with better behaviour for the $p=2$ approximation thanks to the richer approximation space and the presence of interior nodes in the element.

\subsection{Robustness with respect to the interior penalty parameter $\beta$ }

The effect of the interior penalty parameter in the accuracy of the numerical solution is studied next, with the 2D example and meshes of the previous section. Following Remark \ref{rem:beta}, the parameter is taken as \eqref{eq:beta},
with different orders of magnitude for $\alpha$, independent of $h$.
\begin{figure}[!htb]
    \centering
    \includegraphics[width=\textwidth]{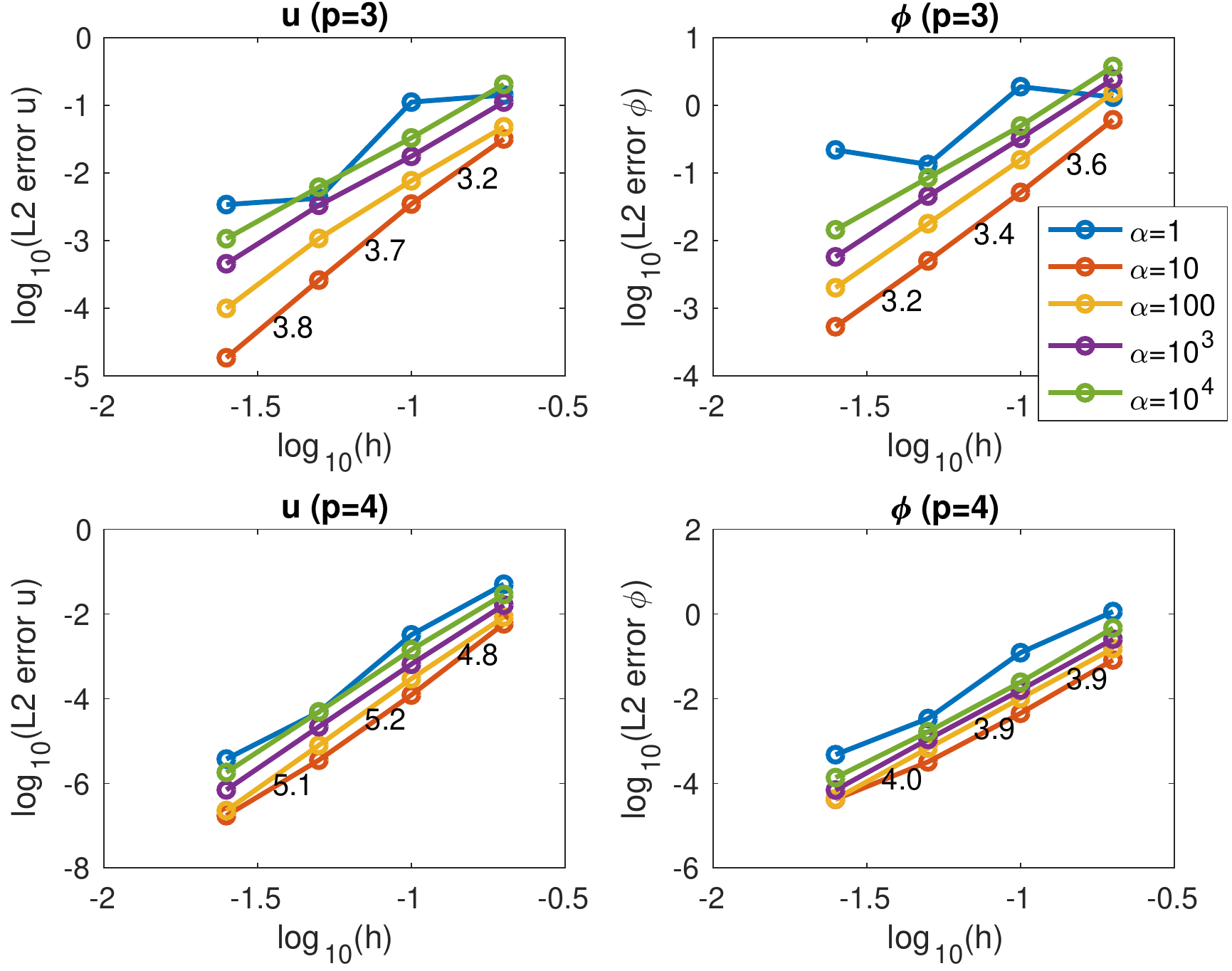}
    \caption{Effect of the $\beta$ parameter in the solution of the \emph{flexoelectricity} problem: convergence plots for the  example in section \ref{sec:2Dconvergence} with $\beta=\alpha El^2/h$ and different values of $\alpha$,  for degree $p=3$ (top) and $p=4$ (bottom).}
    \label{fig:betaDependencyFlexo}
\end{figure}

Figure \ref{fig:betaDependencyFlexo} shows the convergence plots for the flexoelectricity coupled problem, for the displacement $\bu$ (left) and for the potential $\phi$ (right), for degree $p=3$ (top) and $p=4$ (bottom). The slopes of the segments are shown for the plots with $\alpha=10$ for $p=3$, and with $\alpha=100$ for $p=4$. We can observe the poor performance of the method for $\alpha=1$, due to the fact that it is not large enough for a coercive mechanical bilinear form. 

For degree $p=3$, $\alpha=10$ is large enough and provides the best results. Larger values of $\alpha$, several orders of magnitude larger, also lead to high-order convergence, proving the robustness of the method; but, in agreement with the analysis in \cite{Brenner2005} the convergence rate {slowly} decreases for increasing $\alpha$.

Looking to the results for $p=4$ we can observe that, for $\alpha=10$, the bilinear form is coercive for the first meshes, because the elasticity part dominates in the coefficients of the matrix. This is not the case for the last mesh, where higher order terms become more relevant. With $\alpha\geq 100$ the condition in Remark \ref{rem:beta} is satisfied and convergence {is close to} $p+1=5$ for $\bu$ and around $p=4$ for $\phi$, with {almost no} loss in the accuracy for increasing $\beta$.

Thus, from this experiment we conclude that C0-IPM with degree $p=4$ provides excellent results, with convergence rates close to $p+1=5$ for $\bu$ and around $p=4$ for $\phi$,  and {with} little dependency on the particular value of $\beta$, for $\beta\geq 100 El^2/h$.

\begin{figure}[!htb]
    \centering
    \includegraphics[width=\textwidth]{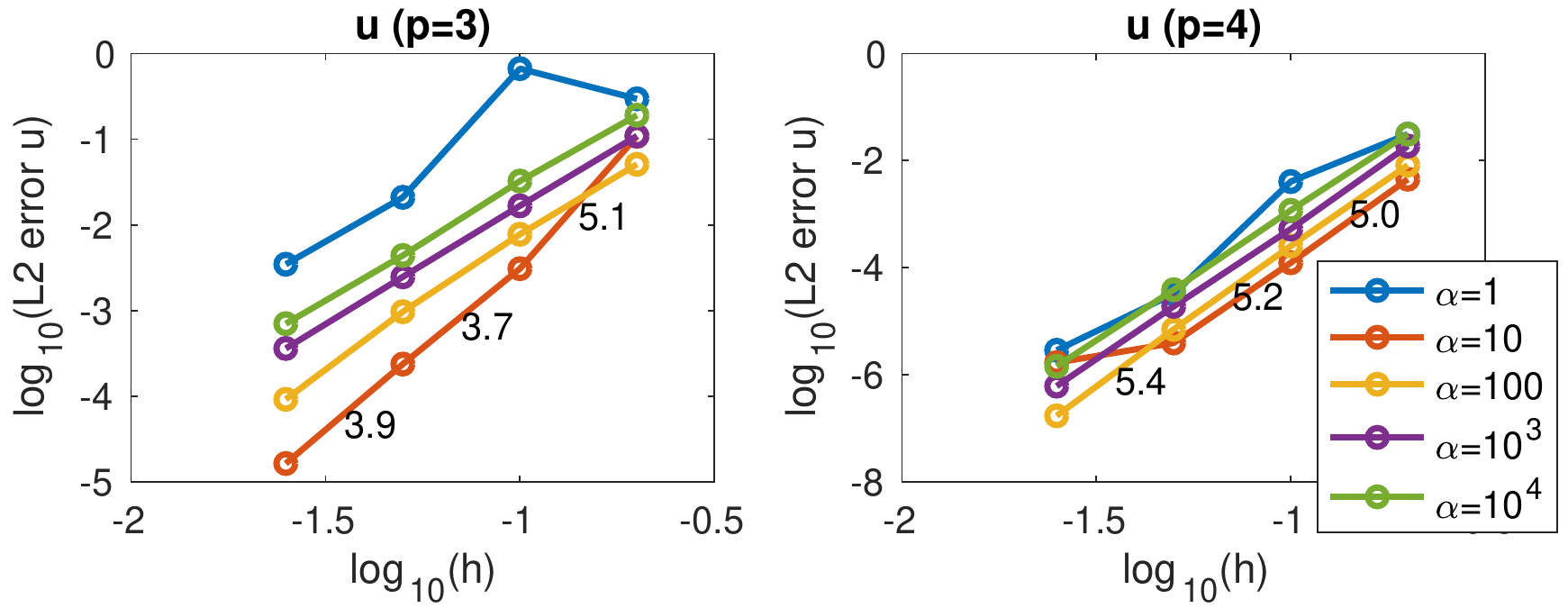}
    \caption{Effect of the $\beta$ parameter in the solution of the \emph{strain gradient elasticity} problem: convergence plots for the displacement, for the  example in section \ref{sec:2Dconvergence}, with $\beta=\alpha El^2/h$ and different values of $\alpha$,  for degree $p=3$ (left) and $p=4$ (right).}
    \label{fig:betaDependencySG}
\end{figure}

The same analysis is performed now for strain gradient elasticity. Figure \ref{fig:betaDependencySG} shows the convergence plots for the displacement $\bu$ for degree $p=3$ (left) and degree $p=4$ (right), with the same conclusions.

\subsection{Cantilever beam}

\begin{figure}[!htb]
    \centering
    \includegraphics[width=0.85\textwidth]{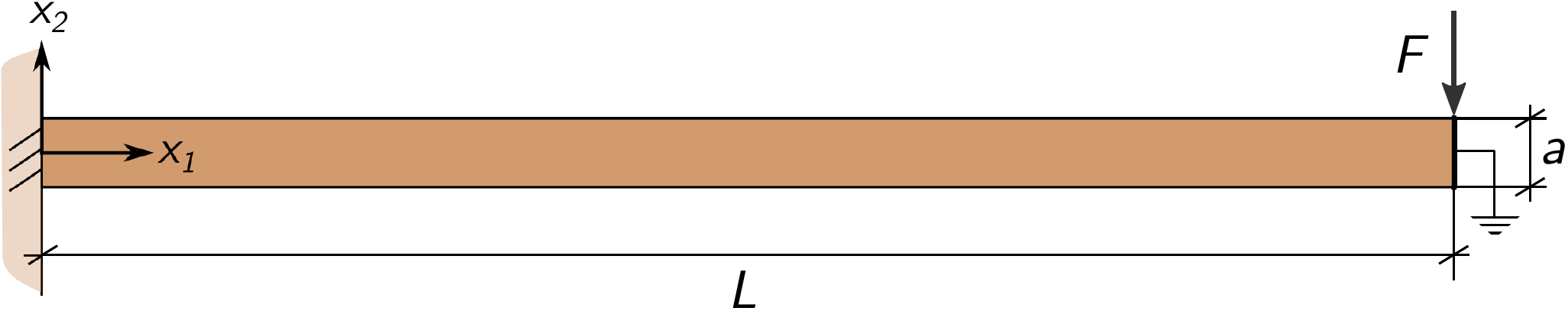}
    \caption{Cantilever beam {under bending and open-circuit boundary conditions}.}
    \label{fig:cantileverBeam}
\end{figure}

The cantilever beam depicted in figure \ref{fig:cantileverBeam} is considered. The aspect ratio is 20, and the width $a$ varies to show the size-dependent nature of flexoelectricity.
The beam is fixed to a wall and grounded on its left end, and it undergoes a punctual force $F$ at the top-right corner. 
The boundary conditions are thus 
\begin{equation} \label{eqn:cantilever}
\begin{array}{rll}
        \bm{u}=&\bm{0}  \qquad & \text{  at  } x_1=0\\
        {j}_2(\bu,\phi)= &-F & \text{  at  } \bm{x}=(L,a/2)\\
        \phi= &0 & \text{  at  } x_1=L,
\end{array}
\end{equation}
where $L=20a$ is the beam length.
To reproduce the results obtained in \cite{Codony2019} with B-splines, the material parameters are
\begin{equation}\label{eq:materialParametersBeam1}
\begin{array}{l}
	 E = 100\, \si{\giga\pascal},\; 
	,\; \kappa_{11} = \kappa_{22} = 11 \, \si{\nano\joule\per\square\volt\per \meter}, \\
	 e_T = -4.4\, \si{\joule\per \volt\per \square\meter},\; \mu_T = 1\, \si{\micro\joule\per \volt\per \meter}, \\
	 l=\nu = \mu_L = \mu_S = e_L = e_S = 0,
\end{array}
\end{equation}
and the piezoelectric principal direction is $x_2$. A uniform discretization with $2\times 2 \times 40$ triangular elements (with characteristic element size $h=0.5a$) of degree $p=4$, and with $\beta=100$, is considered. Since $l=0$, any positive value of $\beta$ provides good results.
\begin{figure}
    \centering
        \includegraphics[width=0.45\textwidth]{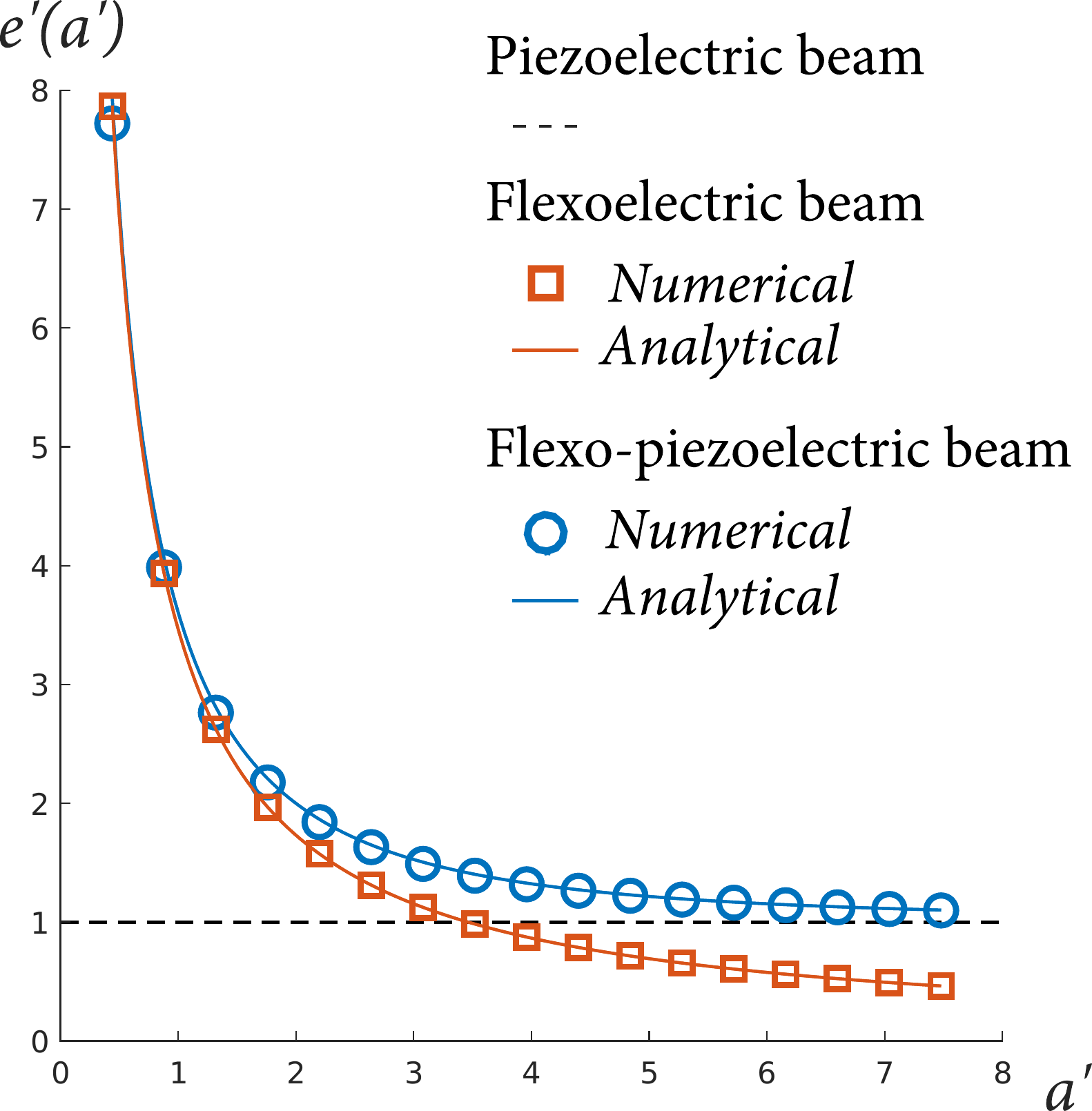}
        \hfill
        \includegraphics[width=0.5\textwidth]{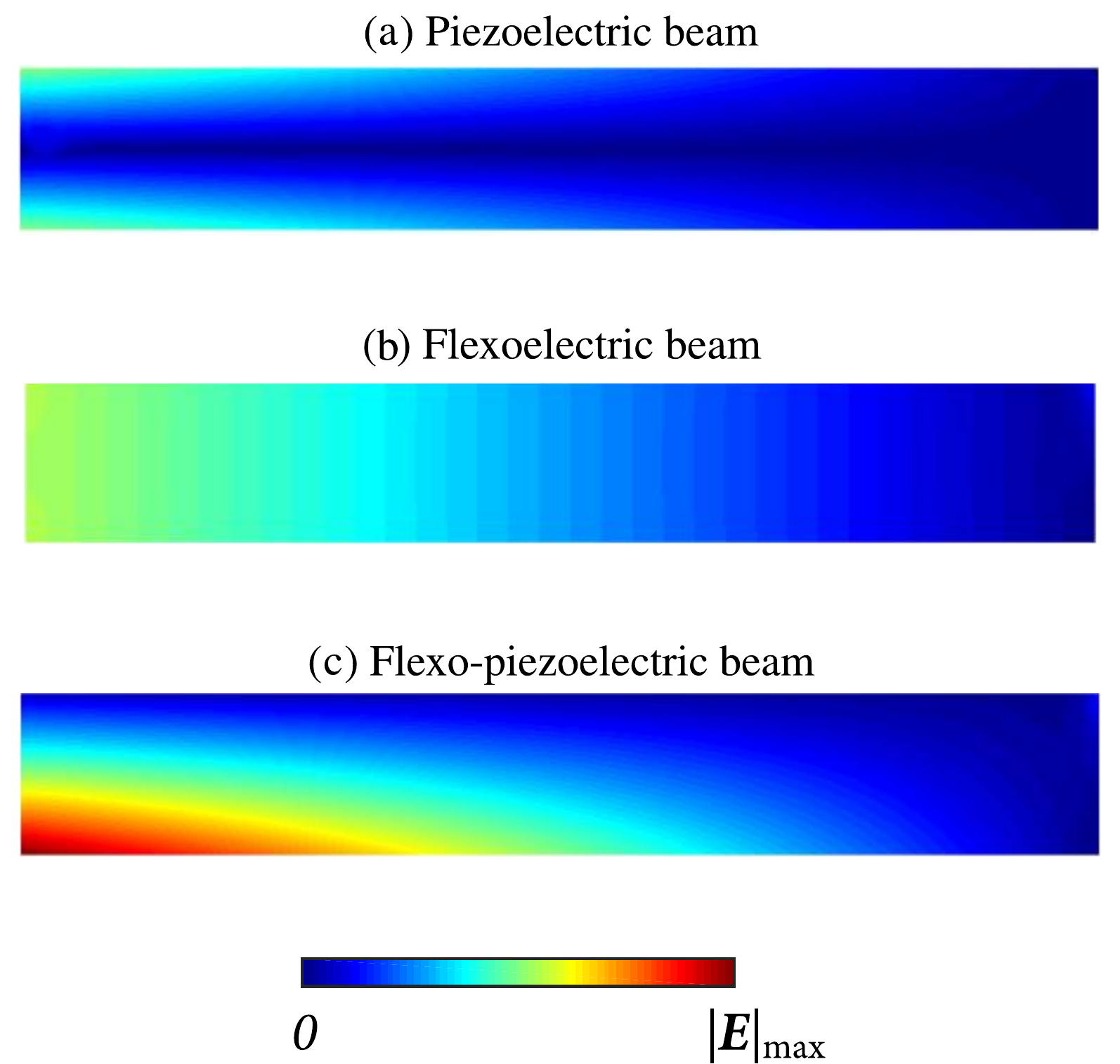}
    \caption{Cantilever beam: (left) normalised effective piezoelectric constant $e'$ as a function of the normalised beam thickness $a'$, and (right) electric field modulus $|\bm{E}|$ with $a'=1.76$ for piezoelectric (a), pure flexoelectric (b) and flexo-piezoelectric (c) beams.}\label{fig:beamResults1}
\end{figure}

Figure \ref{fig:beamResults1} left shows the normalised effective piezoelectric constant{, $e'$,} versus the normalised beam thickness{, $a'$,} defined as
\begin{equation*}
	a'=-ae_T\mu_T^{-1},\quad e':=\cfrac{k_{\text{eff}}}{k_{\text{eff}}|_{\bm{\mu}= \bm{0}}},\quad 	k_{\text{eff}}:=\sqrt{\cfrac{\int_\Omega\bm{E}\cdot\bm{\kappa}\cdot \bm{E}\,\text{d}\Omega}{\int_\Omega\bm{\varepsilon}\cdot\bm{C}\cdot \bm{\varepsilon}\,\text{d}\Omega}},
\end{equation*}
where $k_{\text{eff}}|_{\bm{\mu}= \bm{0}}$ is the effective piezoelectric constant in the absence of flexoelectric effects, i.e. with $\bm{\mu}=\bm{0}$.

The results are in perfect agreement with the B-spline results in \cite{Codony2019}, and with the analytical approximation in \cite{Majdoub2008}:
\begin{equation*}\label{eq:analytical}
	e'|_{\text{flexo}}(a')\simeq\sqrt{\cfrac{12}{{a'} ^2}},\quad e'|_{\text{flexo-piezo}}(a')\simeq\sqrt{1+\cfrac{12}{{a'}^2}}.
\end{equation*}
The plots in figure \ref{fig:beamResults1}  also illustrate how flexoelectricity is a size dependent phenomenon, with relevant and even crucial effect for very small scales.

\begin{figure}[!b]
    \centering
    \includegraphics[width=0.85\textwidth]{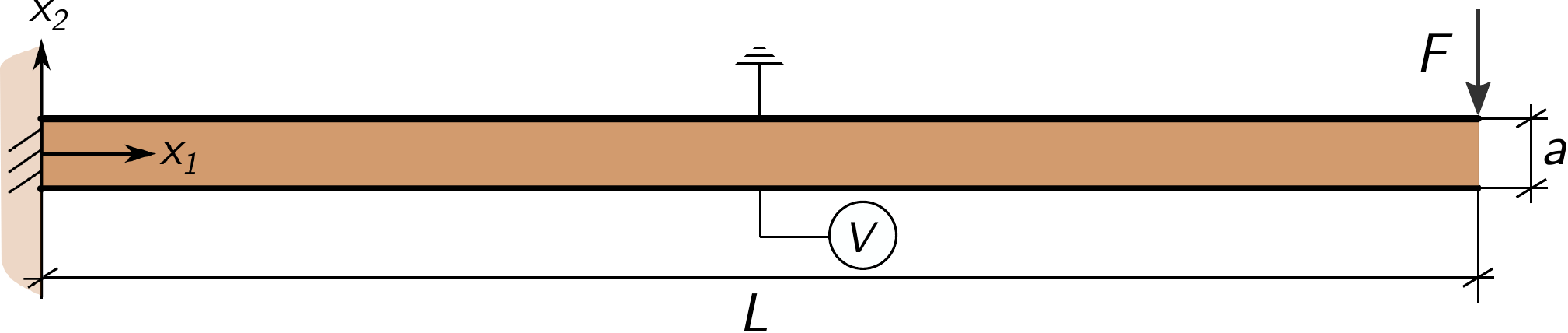}
    \caption{Cantilever beam {under closed-circuit boundary conditions.}}
    \label{fig:cantileverBeamCC}
\end{figure}

\begin{figure}[!b]
    \centering
        \includegraphics[width=0.8\textwidth]{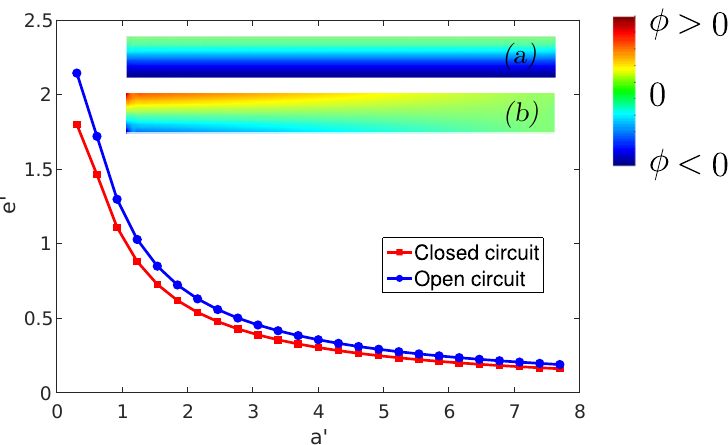}
     \caption{Normalised effective piezoelectric constant as a function of the normalised beam thickness, and example of the distribution of electric potential $\phi$ in a {flexo-piezoelectric} beam with closed (a) and open (b) circuit. Note that the aspect ratio of the beams has been modified to better observe the potential distribution along the beam. }\label{fig:amir}
\end{figure}

\subsection{Open and Closed circuit in the cantilever beam}
For further validation of the C0-IPM computational model, we now consider the open and closed circuit example in \cite{Abdollahi2014}, where maximum-entropy approximations (LME) were used. The problem is solved on the same beam with the same FE mesh.
The material parameters are  now
\begin{equation*}
\begin{array}{l}
	 E = 100\, \si{\giga\pascal},\; \nu=0.37,\\
 \kappa_{11} = 11 \si{\nano\joule\per\square\volt\per \meter},\; \kappa_{22} = 12.48\; \si{\nano\joule\per\square\volt\per \meter}, \\
	 e_T =-4.4\, \si{\joule\per \volt\per \square\meter},\; \mu_T = \mu_L = 1\, \si{\micro\joule\per \volt\per \meter}, \\
	 l=  \mu_S = e_S =  e_L= 0,
\end{array}
\end{equation*}
and the the piezoelectric principal direction is again $x_2$.
The mechanical boundary conditions are the same as in the previous case.

For the electrostatic boundary conditions, two different cases are considered: open and closed circuit. The open circuit is the one considered in the previous example, with grounded right end, that is $\phi=0$ at $x_1=L${, as shown in figure \ref{fig:cantileverBeam}}. In the closed circuit, the upper side is grounded and an electrode is placed on the bottom side, that is
\begin{equation*}
\phi=\left\lbrace
\begin{array}{ll}
0 &\text{ for }x_2=a/2 \\
V &\text{ for } x_2=-a/2 ,
\end{array}
\right.
\end{equation*}
where $V$ is a free constant value, see figure \ref{fig:cantileverBeamCC}. The electrode condition is enforced setting all potential nodal values on the bottom boundary to be equal to the first one, {with Lagrange multipliers in our implementation.} 

Figure \ref{fig:amir} shows the normalised effective piezoelectric constant $e'$ as a function of the normalised thickness $a'$. Again, we observe that flexoelectricy becomes relevant for small scales.  For the open circuit, comparing to the previous results in figure \ref{fig:beamResults1}, where $\nu=\mu_L = 0$, this more general model gives lower values for the normalised effective piezoelectric constant. On other hand, the open circuit setting leads to larger values of the effective piezoelectric constant. The numerical results are  in perfect agreement with the LME results in \cite{Abdollahi2014} demonstrating again the applicability of C0-IPM for the study and design of flexoelectric devices.

\subsection{Actuator example}

In this section we consider an actuator beam also from \cite{Abdollahi2014}. The displacement is fixed on the left boundary, and a potential difference is applied at the top and bottom sides, leading to a bending of the beam. That is,
\begin{equation*}
\bm{u}(0,x_2)=\bm{0},\quad
 \phi \left( x_1,\frac{a}{2} \right)=0, \quad
 \phi\left( x_1,-\frac{a}{2} \right)=V,
\end{equation*}
on the same beam, i.e. $ \Omega=(0,20a)\times(-a/2,a/2)$. 
\begin{figure}[b!]
    \centering
    \includegraphics[width=0.75\textwidth]{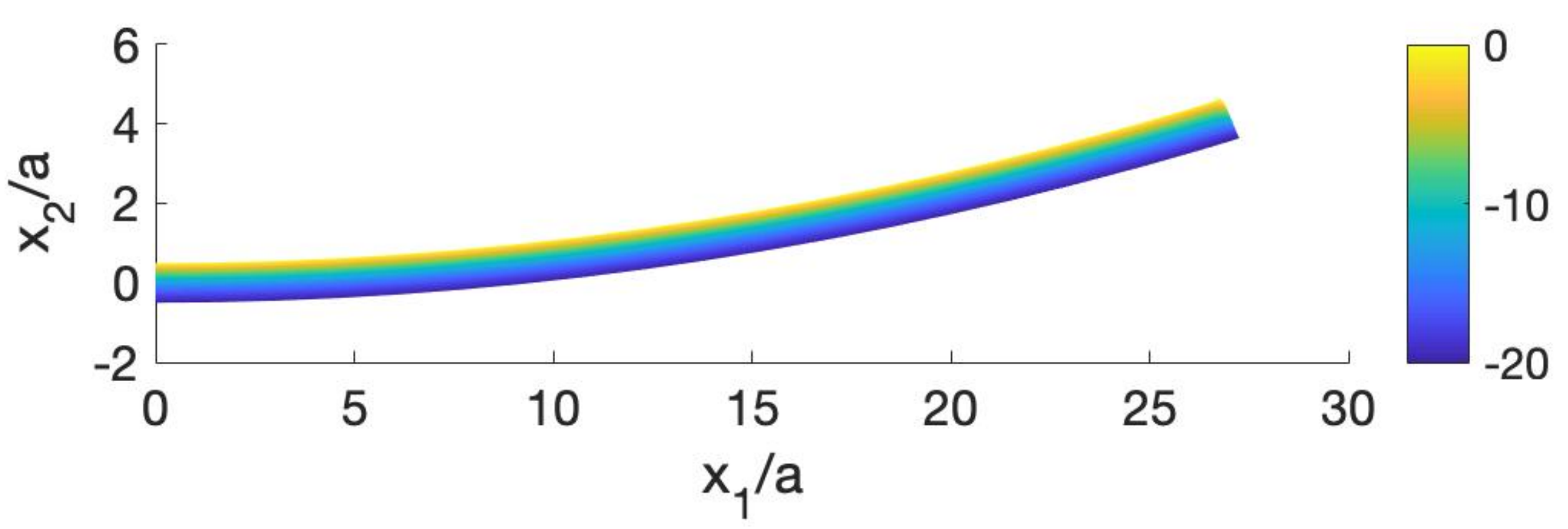}
    \caption{Actuator beam {with width $a=2.5\,\si{\micro \meter}$ and} $l=0$: deformed beam and potential.}
    \label{fig:actuator_deformed}
\end{figure}
The material parameters are \eqref{eq:materialParametersBeam1} and the applied voltage is $V=-8a \;\si{\mega \volt}$.

\begin{figure}[p]
    \centering
    \includegraphics[width=0.75\textwidth]{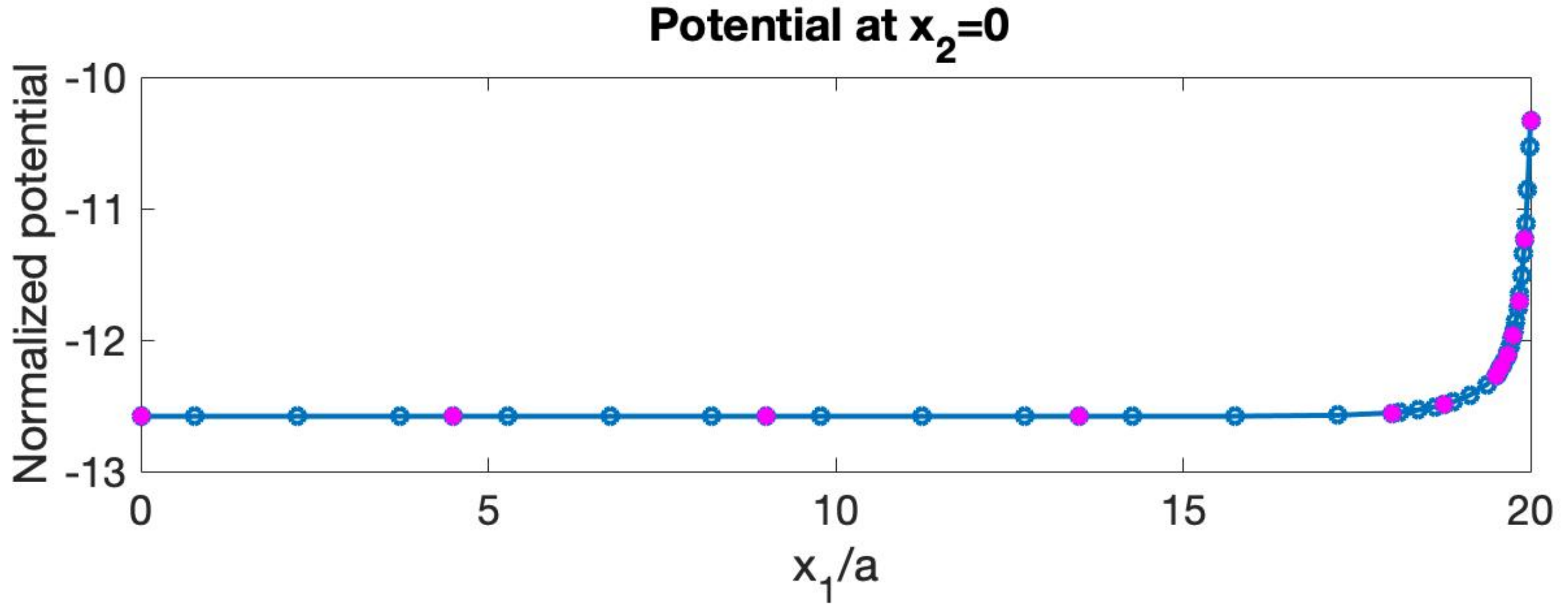}
    \caption{Actuator beam {with width $a=2.5\,\si{\micro \meter}$ and} $l=0$: normalised potential ($\phi/250)$ along the $x_2=0$ horizontal mid section. The potential presents a sharp variation close to the right boundary. Magenta dots correspond to the boundary of the $p=4$ quadrilaterals.}
    \label{fig:actuator_potentialy0}
\end{figure}

\begin{figure}[p]
    \centering
    \includegraphics[width=0.45\textwidth]{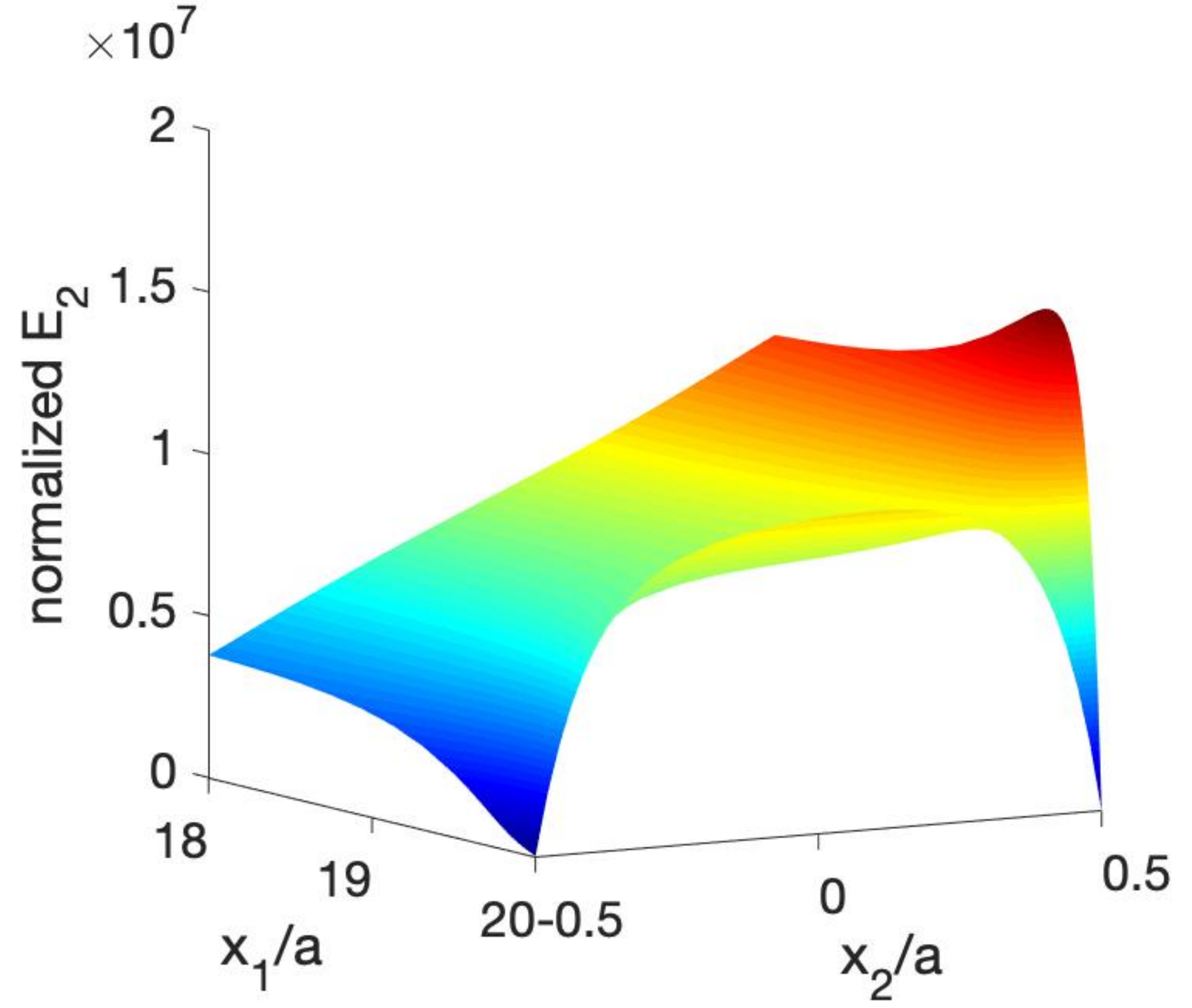} \hfill
    \includegraphics[width=0.5\textwidth]{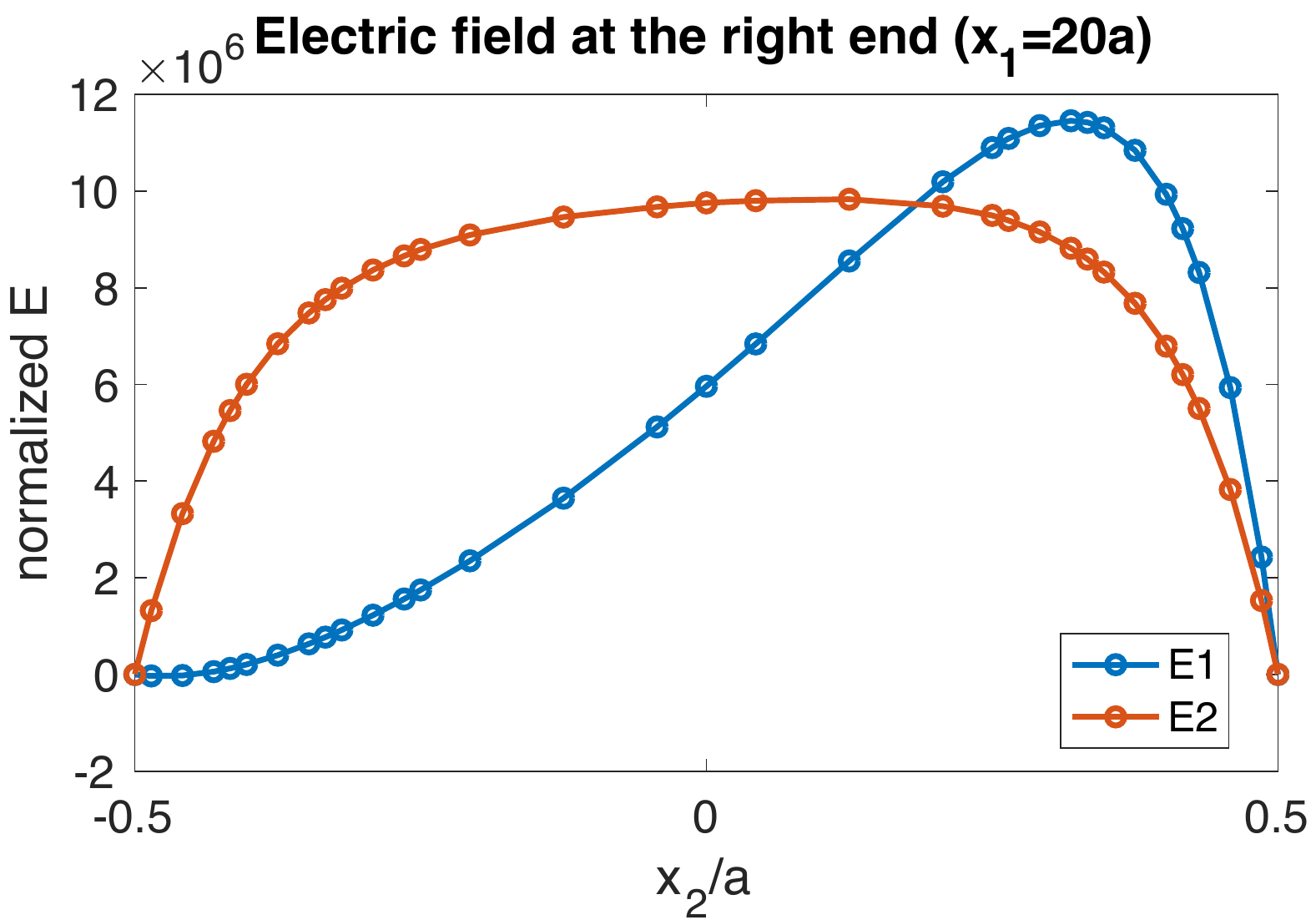}
    \caption{Actuator beam with $l=0$: detail of the normalised vertical electric field ($E_2/10^9$)  close to the right end, and plot of the components of the normalised electric field at the right end ($x_1=20a$).}
    \label{fig:actuator_E2detail}
\end{figure}

Figure \ref{fig:actuator_deformed} shows the potential on the deformed beam  for width $a=2.5\,\si{\micro \meter}$. The potential seems to be smooth, but the section along $x_2=0$ in figure \ref{fig:actuator_potentialy0} reveals a sharp variation close to the right end. Consequently, the electric field also presents sharp variations close to the right end, as shown in figure \ref{fig:actuator_E2detail}. 
\begin{figure}[p]
\centering
    \centering
    \includegraphics[width=\textwidth]{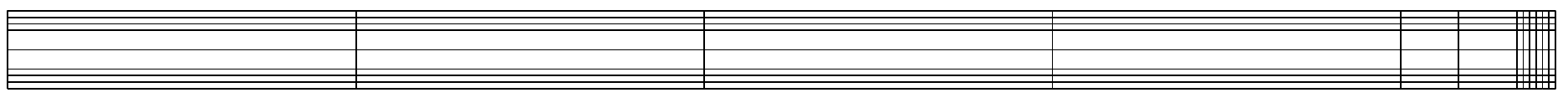}
    \includegraphics[width=\textwidth]{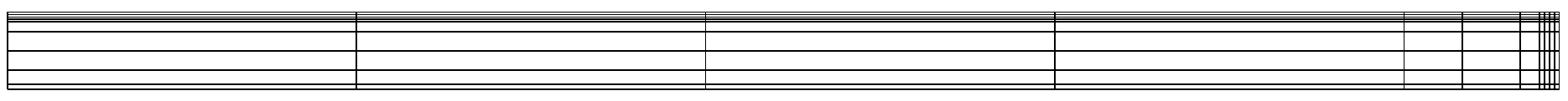}
    \caption{Adapted p=4 quadrilateral FE mesh for $l=0$ with $\min(h)=a/12$ (top) and for $l=0.1a$ with $\min(h)=a/2^5$ (bottom).}
    \label{fig:actuator_mesh}
\end{figure}
\begin{figure}[!htb]
    \centering
\vspace{0.5cm}
    \includegraphics[width=0.45\textwidth]{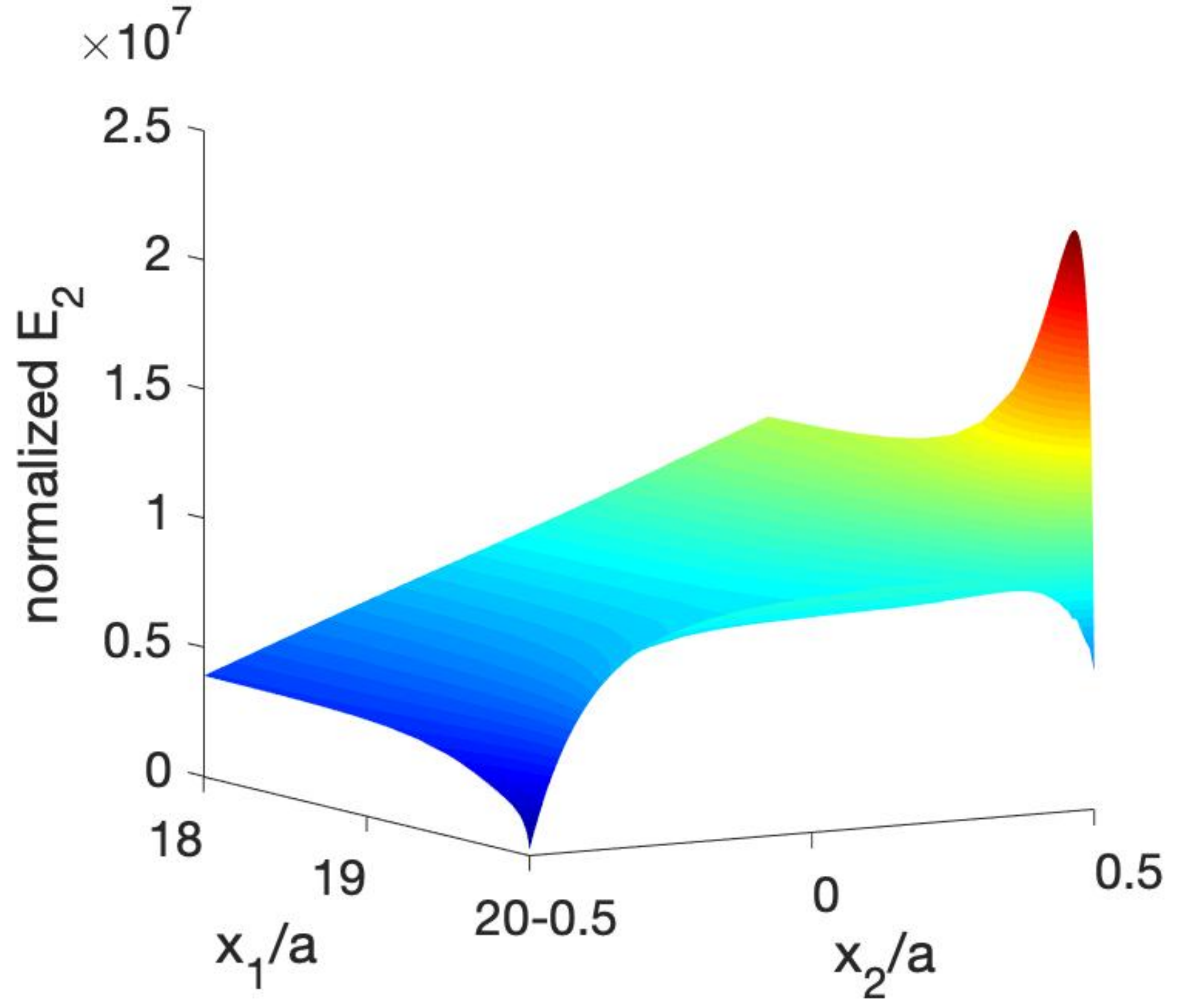} \hfill
    \includegraphics[width=0.5\textwidth]{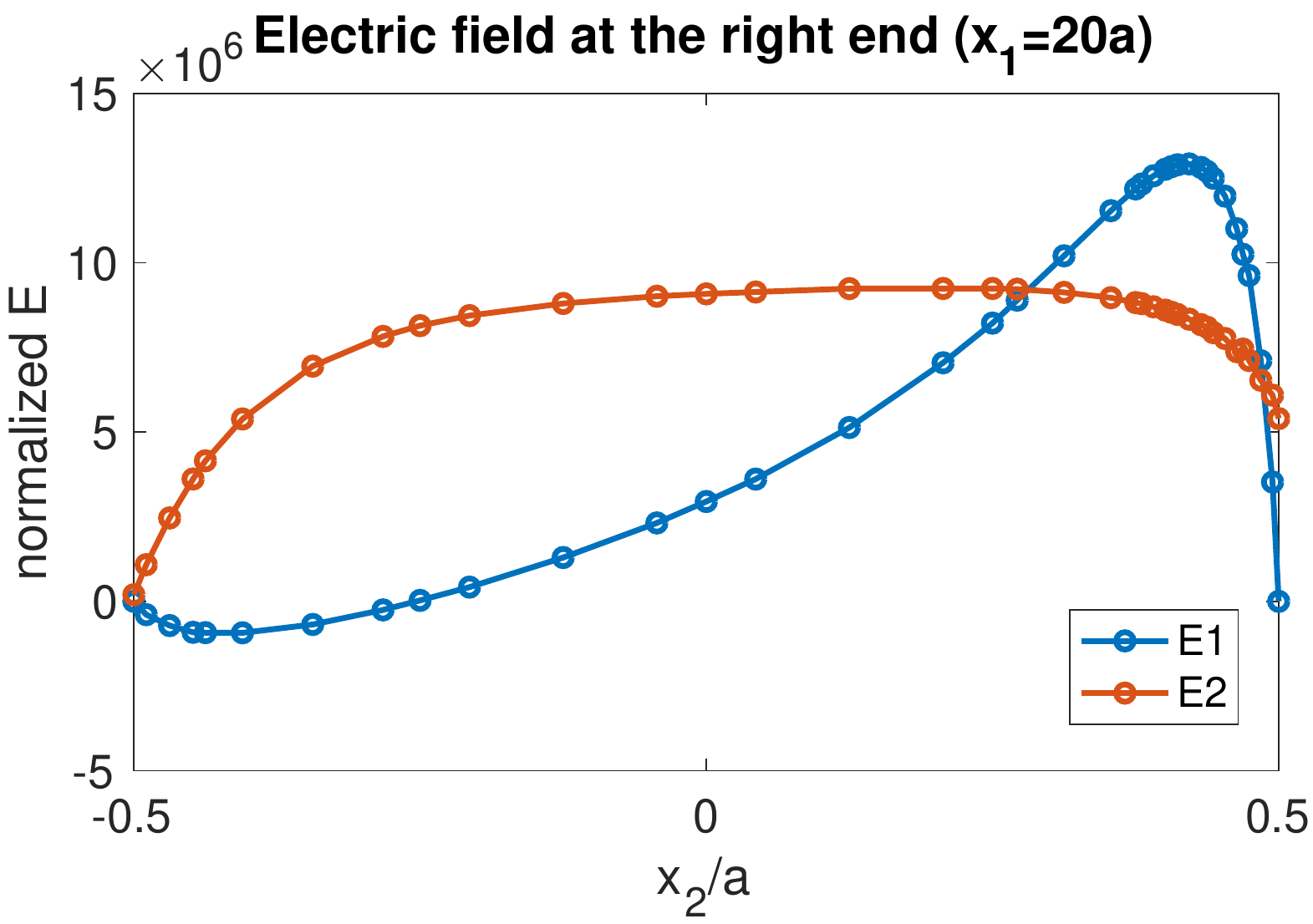}
    \caption{Actuator beam with $l=0.1a$: detail of the normalised vertical electric field $E_2$  close to the right end, and plot of the components of the normalised electric field at the right end ($x_1=20a$).}
    \label{fig:actuator_E2detail_l1}
\end{figure}

These results have  been computed on the adapted quadrilateral mesh in figure \ref{fig:actuator_mesh} (top), with {degree $p=4$} and $\beta=1$. The mesh has been refined to capture the sharp variations in the solution; otherwise, numerical oscillations spoil the solution in the whole domain. It is also worth mentioning that the plot in figure \ref{fig:actuator_E2detail} coincides in magnitude and shape with the results in \cite{Abdollahi2014} with LME, but getting rid of the smooth oscillations.

Sharp variations along the boundary in the solution of flexoelectricity problems can be even more pronounced, as can be observed in figure \ref{fig:actuator_E2detail_l1}. In this case the problem is solved with strain gradient elasticity, with $l=0.1a$,  on the $p=4$ adapted mesh in figure \ref{fig:actuator_mesh} with $\min(h)=a/2^5$; with smaller element size along the boundary to capture the high curvatures in the electric field. The stabilization parameter is again taken as $\beta=100El^2/\min(h)$, providing stable results.

\subsection{Periodicity}

The implementation of periodicity boundary conditions in the C0-IPM method is straightforward, by simply considering the periodicity faces as interior faces and imposing the periodicity constraints on the boundary nodal values. 
Considering the periodicity faces as interior faces, that is in $\I$, ensures that $\mathcal{C}^1$ continuity is enforced in weak form and that internal forces are equilibrated also on the periodicity boundary. The periodicity conditions for the nodal values can be implemented, for instance, by means of Lagrange multipliers, or reducing the system to the periodic space.
\begin{figure}[!htb]
    \centering
    \includegraphics[width=0.9\textwidth]{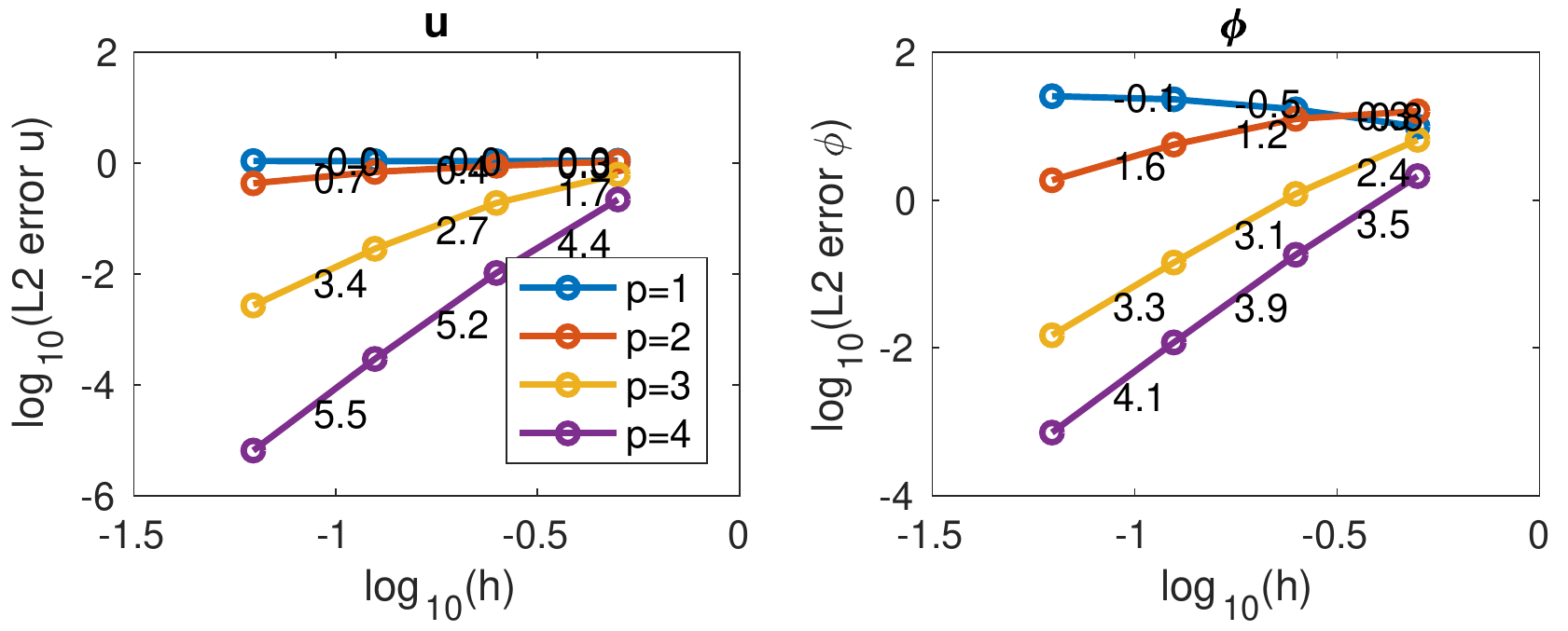} \hfill
    \caption{Convergence test for the solution of the flexoelectricity equations in a square domain, with periodicity in the $x_1$ direction.}
    \label{fig:convergence_periodic}
\end{figure}

As a verification example, figure \ref{fig:convergence_periodic} shows the evolution of the error under nested refinement for the solution of the flexoelectricity coupled problem \eqref{eq:problemStatement} in a square domain $\Omega=(0,1)^2$ with a regular triangular mesh. First Dirichlet and second Neumann conditions, \eqref{eq:Dirichlet1} and \eqref{eq:Neumann1},  are set on the top and bottom boundaries, and periodicity is imposed in the $x_1$ direction. That is, \eqref{eq:interfaceConditions} is imposed identifying the left and right boundary as the same boundary and including it in $\I$. {The body force $\bm{b}$, the free charge $q$, and the data for the boundary conditions on the top and bottom boundaries, are set so that the analytical solution is \eqref{eq:sinAnalytical}. The stabilization parameter is \eqref{eq:beta} with $\alpha=100$.}

The errors exhibit the same behavior as {in the convergence analysis in section \ref{sec:2Dconvergence}.}

\subsection{3D convergence test}

\begin{figure}[!p]
    \centering
    \includegraphics[width=0.3\textwidth]{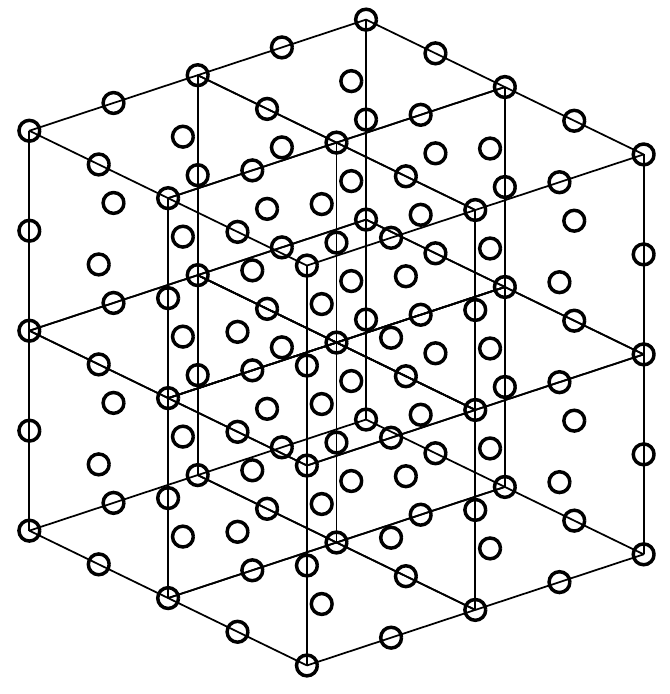} \hfill
    \caption{3D mesh for degree $p=2$ and second level of refinement.}
    \label{fig:mesh3d}
\end{figure}

The flexoelectricity equations are now solved in a cube, $\Omega=(0,0.5)^3$, to show the applicability of the method also in 3D. The mesh for degree $p=2$ and the second level of refinement is shown in figure \ref{fig:mesh3d}. First Dirichlet and second Neumann boundary conditions are considered in the whole boundary, and the material parameters are \eqref{eq:materialParameters}. The data is set so that the solution is 
\begin{align*}
&\bu=\left[ \cos(2 \pi(x_1+2x_2-x_3)),\sin(2 \pi(x_1+2x_2-x_3)),\cos(2 \pi(x_1+2x_2-x_3))\right]^T, \\ 
&\phi=\sin(2 \pi(x_1+2x_2-x_3)).
\end{align*}
\begin{figure}[!p]
    \centering
    \includegraphics[width=0.82\textwidth]{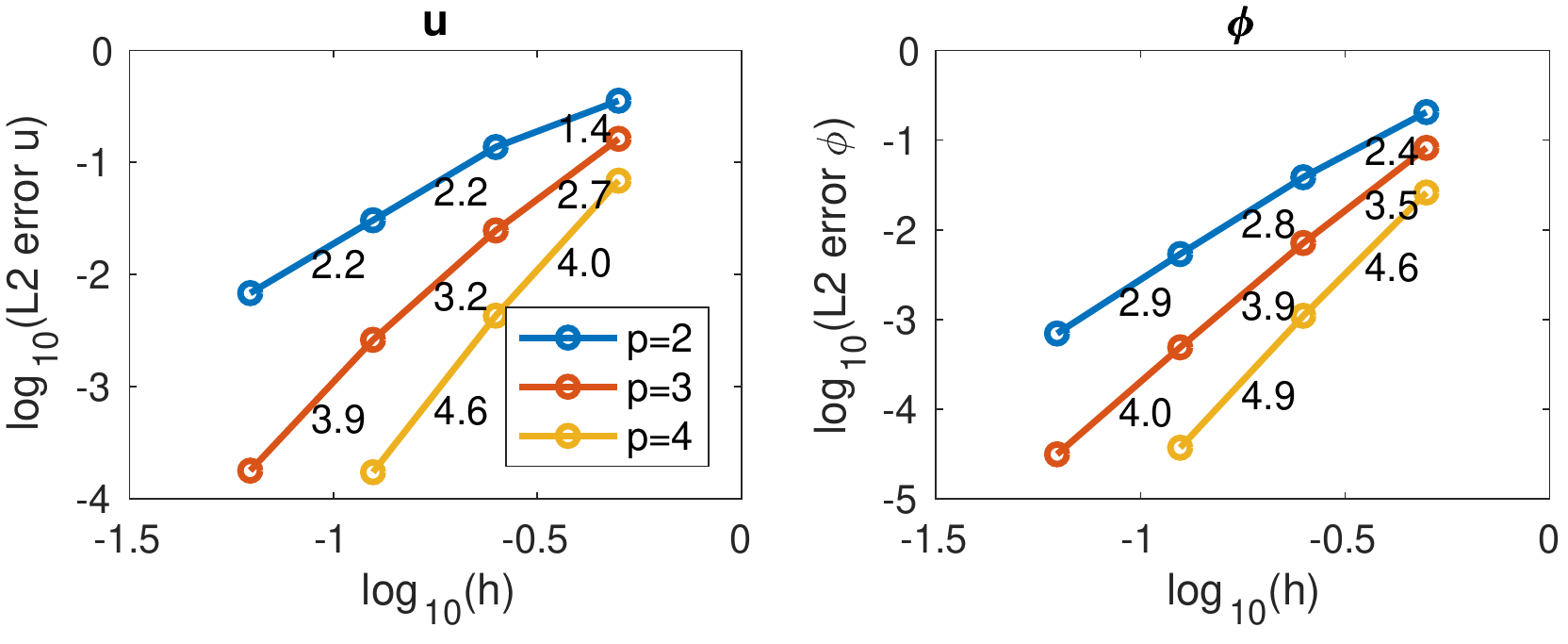} \hfill
    \caption{3D solution in a cube: convergence plots for \emph{strain gradient elasticity} and electric potential (decoupled problem).}
    \label{fig:errors_3d_alpha100_uncoupled}
\end{figure}
\begin{figure}[!p]
    \centering
       \includegraphics[width=0.82\textwidth]{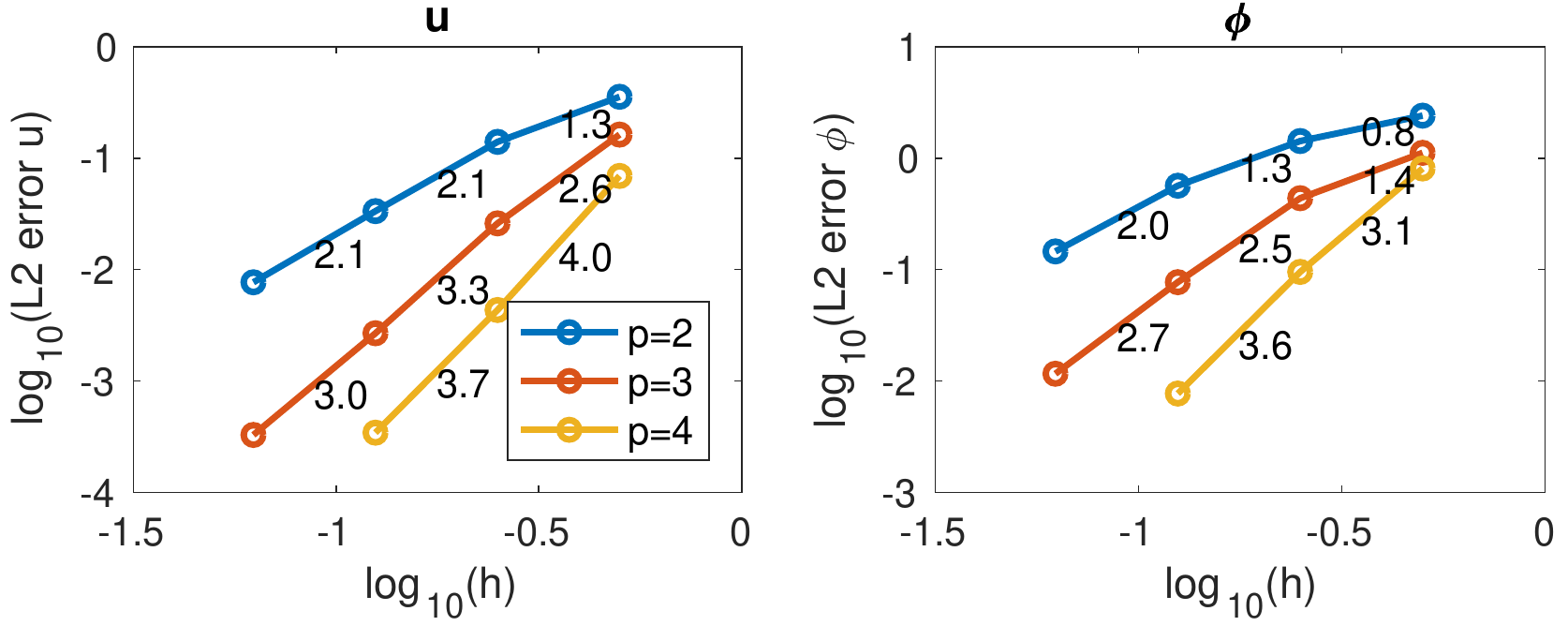} \hfill
    \caption{3D solution in a cube: convergence plots for the coupled \emph{flexoelectricity} problem.}
    \label{fig:errors_3d_alpha100}
\end{figure}

Figures \ref{fig:errors_3d_alpha100_uncoupled} and \ref{fig:errors_3d_alpha100} show the convergence plots with $\beta=100El^2/h$, for strain gradient elasticity (solving the decoupled problem) and for flexoelectricity, respectively. As in 2D, the method does not converge for degree $p=1$; thus, we show the results for $p=2,3,4$. 

Robust high-order convergence is observed in all cases, providing accurate results. Again, in agreement with the analysis in \cite{Brenner2005}, the convergence is suboptimal; but still with order close to $p+1$ for the displacement in the strain gradient elasticity problem for $p\geq 3$. Again, we also observe that the flexoelectricy coupling provokes a loss in the convergence rate and the accuracy of the solution; with order close to $p$ for the displacement and the potential in this example.

Further numerical experiments show that, with these regular hexahedra meshes, the error has very little dependency on the particular value of $\beta \geq 100El^2/h$. 

\section{Conclusions}
A novel C0-IPM formulation for strain gradient elasticity and flexoelectricity is proposed. The weak form involves second derivatives of the displacement in the interior of the elements, plus integrals on the element faces, weakly imposing continuity of the displacement derivatives, as well as equilibrium of internal forces across element faces and on interior edges (vertexes in 2D). 

The formulation is stable, with a symmetric and positive definite matrix for the strain gradient elasticity operator, for large enough interior parameter $\beta$. An eigenvalue problem is stated to determine a bound for $\beta$, which leads to a general formula for the parameter: $\beta=\alpha El^2/h$, with constant $\alpha$ independent of the element size. Thus, differently to non-consistent penalty methods, and as usual in interior penalty methods, moderate values for $\beta$ provide stable and accurate results.

Standard $\mathcal{C}^0$ FE approximations are considered, retaining the advantages and computational efficiency of high-order FE. The implementation is based on  assembly of elemental matrices, with standard FE numerical integration and nodal approximation, the discretization can be adapted to the geometry and locally refined where needed, no additional unknowns are needed, and material interfaces can be directly considered just adapting the mesh, as usual in FE computations.  

The application of C0-IPM to problems with periodicity boundary conditions is straightforward, just considering the periodicity faces as interior faces (thus, imposing $\mathcal{C}^1$ continuity and equilibrium of forces in weak form) and setting the periodicity conditions on the nodal values.

Convergence tests, on 2D non-uniform curved triangular meshes and on 3D hexahedra regular meshes, show high-order convergence of the method for degree $p\geq 3$. A slow continuous loss in the convergence rate for increasing $\beta$ is observed for $p=3$, which is in agreement with the analysis for the biharmonic equation in \cite{Brenner2005} and the results for Kirchhoff plates in \cite{Fojo2020}. Fortunately, for $p=4$ the convergence shows little dependency on $\beta$. In any case, in all examples, the convergence rates are at least close to $p$ for both variables,  demonstrating the good behaviour of the method for $p\geq 3$.

The computational tool is also validated by comparison with previous works solving realistic actuator and sensor problems on a beam, with perfect agreement.

\section*{Acknowledgements}
This work was supported by the European Research Council (StG-679451 to Irene Arias) and Generalitat de Catalunya (2017-SGR-1278)

\bibliographystyle{unsrt}
\bibliography{BibliografiaFlexo}

\begin{thebibliography}{10}

\bibitem{Shu2019}
Longlong Shu, Renhong Liang, Zhenggang Rao, Linfeng Fei, Shanming Ke, and
  Yu~Wang.
\newblock Flexoelectric materials and their related applications: A focused
  review.
\newblock {\em Journal of Advanced Ceramics}, 8(2):153--173, 2019.

\bibitem{Mao2016}
Sheng Mao, Prashant Purohit, and N.~Aravas.
\newblock Mixed finite-element formulations in piezoelectricity and
  flexoelectricity.
\newblock {\em Proceedings of the Royal Society A: Mathematical, Physical and
  Engineering Science}, 472:20150879, 2016.

\bibitem{Deng2017}
Feng Deng, Qian Deng, Wenshan Yu, and Shengping Shen.
\newblock {Mixed Finite Elements for Flexoelectric Solids}.
\newblock {\em Journal of Applied Mechanics}, 84(8), 2017.

\bibitem{Abdollahi2014}
Amir Abdollahi, Christian Peco, Daniel Mill{\'a}n, Marino Arroyo, and Irene
  Arias.
\newblock Computational evaluation of the flexoelectric effect in dielectric
  solids.
\newblock {\em Journal of Applied Physics}, 116:093502--093502, 09 2014.

\bibitem{Ghasemi2017}
Hamid Ghasemi, Harold~S. Park, and Timon Rabczuk.
\newblock A level-set based {IGA} formulation for topology optimization of
  flexoelectric materials.
\newblock {\em Computer Methods in Applied Mechanics and Engineering}, 313:239
  -- 258, 2017.

\bibitem{Nanthakumar2017}
S.S. Nanthakumar, Xiaoying Zhuang, Harold~S. Park, and Timon Rabczuk.
\newblock Topology optimization of flexoelectric structures.
\newblock {\em Journal of the Mechanics and Physics of Solids}, 105:217 -- 234,
  2017.

\bibitem{Sevilla2011}
R.~Sevilla and S.~Fern{\'a}ndez-M{\'e}ndez.
\newblock Numerical integration over {2D} {NURBS}-shaped domains with
  applications to {NURBS}-enhanced \{FEM\}.
\newblock {\em Finite Elements in Analysis and Design}, 47(10):1209 -- 1220,
  2011.

\bibitem{Codony2019}
D.~Codony, O.~Marco, S.~Fern{\'a}ndez-M{\'e}ndez, and I.~Arias.
\newblock An immersed boundary hierarchical {B}-spline method for
  flexoelectricity.
\newblock {\em Computer Methods in Applied Mechanics and Engineering}, 354:750
  -- 782, 2019.

\bibitem{dePrenter2017}
F.~{de Prenter}, C.V. Verhoosel, G.J. {van Zwieten}, and E.H. {van Brummelen}.
\newblock Condition number analysis and preconditioning of the finite cell
  method.
\newblock {\em Computer Methods in Applied Mechanics and Engineering}, 316:297
  -- 327, 2017.
\newblock Special Issue on Isogeometric Analysis: Progress and Challenges.

\bibitem{Arnold1982}
Douglas~N. Arnold.
\newblock An interior penalty finite element method with discontinuous
  elements.
\newblock {\em SIAM Journal on Numerical Analysis}, 19(4):742--760, 1982.

\bibitem{Nitsche1971}
J.~Nitsche.
\newblock Über ein variationsprinzip zur lösung von dirichlet-problemen bei
  verwendung von teilräumen, die keinen randbedingungen unterworfen sind.
\newblock {\em Abhandlungen aus dem Mathematischen Seminar der Universität
  Hamburg}, 36(1):9--15, 1971.

\bibitem{Engel2002}
Gerald Engel, Krishna Garikipati, Thomas Hughes, Mats Larson, Luca Mazzei, and
  R.~Taylor.
\newblock Continuous/discontinuous finite element approximations of
  fourth-order elliptic problems in structural and continuum mechanics with
  applications to thin beams and plates.
\newblock {\em Computer Methods in Applied Mechanics and Engineering},
  191:3669--3750, 2002.

\bibitem{Brenner2005}
Susanne~C. Brenner and Li-Yeng Sung.
\newblock {C}0 {I}nterior {P}enalty {M}ethods for fourth order elliptic
  boundary value problems on polygonal domains.
\newblock {\em Journal of Scientific Computing}, 22(1):83--118, 2005.

\bibitem{Fojo2020}
Dani Fojo, David Codony, and Sonia Fern{\'a}ndez-M{\'e}ndez.
\newblock A {C}0 {I}nterior {P}enalty {M}ethod for 4th order {PDE}s.
\newblock {\em Reports@SCM}, 2020.

\bibitem{Wells2003}
Garth Wells, Krishna Garikipati, and Luisa Molari.
\newblock A discontinuous galerkin method for strain gradient-dependent damage.
\newblock {\em Computer Methods in Applied Mechanics and Engineering},
  193:3633--3645, 2003.

\bibitem{Brenner2012}
Susanne~C. Brenner, Shiyuan Gu, Thirupathi Gudi, and Li-Yeng Sung.
\newblock {A} quadratic {C}0 interior penalty method for linear fourth order
  boundary value problems with boundary conditions of the {C}ahn-{H}illiard
  type.
\newblock {\em SIAM Journal on Numerical Analysis}, 50(4):2088--2110, 2012.

\bibitem{Fernandez2004}
Sonia Fern{\'a}ndez-M{\'e}ndez and Antonio Huerta.
\newblock Imposing essential boundary conditions in mesh-free methods.
\newblock {\em Computer Methods in Applied Mechanics and Engineering},
  193(12):1257 -- 1275, 2004.

\bibitem{Griebel2002}
Michael Griebel and Marc Schweitzer.
\newblock A particle-partition of unity method - part v: Boundary conditions.
\newblock {\em Geometric Analysis and Nonlinear Partial Differential
  Equations}, 41:519--542, 05 2002.

\bibitem{Annavarapu2012}
Chandrasekhar Annavarapu, Martin Hautefeuille, and John~E. Dolbow.
\newblock A robust {N}itsche's formulation for interface problems.
\newblock {\em Computer Methods in Applied Mechanics and Engineering},
  225-228:44 -- 54, 2012.

\bibitem{Chen1995}
Qi~Chen and Ivo Babu\u{s}ka.
\newblock Approximate optimal points for polynomial interpolation of real
  functions in an interval and in a triangle.
\newblock {\em Computer Methods in Applied Mechanics and Engineering},
  128(3):405 -- 417, 1995.

\bibitem{RuizGirones2019}
Eloi Ruiz-Giron\'es, Abel Gargallo-Peir\'o, Josep Sarrate, and Xevi Roca.
\newblock Automatically imposing incremental boundary displacements for valid
  mesh morphing and curving.
\newblock {\em Computer-Aided Design}, 112:47 -- 62, 2019.

\bibitem{Majdoub2008}
M.~Majdoub, Pradeep Sharma, and Tcagin Cagin.
\newblock Enhanced size-dependent piezoelectricity and elasticity in
  nanostructures due to the flexoelectric effect.
\newblock {\em Physical Review B}, 77, 03 2008.

\end{thebibliography}

\end{document}